\documentclass[12pt]{amsart}

\usepackage{hyperref}

\usepackage{geometry}
\usepackage{marginnote}

\usepackage{amssymb}

\usepackage{xcolor}

\DeclareFontFamily{U}{mathx}{\hyphenchar\font45}
\DeclareFontShape{U}{mathx}{m}{n}{<5> <6> <7> <8> <9> <10> <10.95> <12> <14.4> <17.28> <20.74> <24.88> mathx10}{}
\DeclareSymbolFont{mathx}{U}{mathx}{m}{n}
\DeclareFontSubstitution{U}{mathx}{m}{n}
\DeclareMathAccent{\widecheck}{0}{mathx}{"71}
\DeclareMathAccent{\wideparen}{0}{mathx}{"75}

\newcommand\la{\langle}
\newcommand\ra{\rangle}
\newcommand\ov{\overline}

\newcommand\wh{\widehat}
\newcommand\mc{\mathcal}

\newcommand\ganz{\mathbb Z} 
\newcommand\real{\mathbb R}

\newcommand\conv{\text{\rm Conv}}

\newcommand\gen{\text{\rm Span}}
\newcommand\stab{\text{\rm Stab}}

\newcommand\supp{\text{\rm Supp}}

\newcommand\al{\alpha}
\newcommand\be{\beta}
\newcommand\ga{\gamma}
\newcommand\ve{\varepsilon}

\newcommand\smeno{\smallsetminus}

\newtheorem*{thm*}{Theorem}
\newtheorem{thm}{Theorem}[section]
\newtheorem{pro}[thm]{Proposition}
\newtheorem{lem}[thm]{Lemma}
\newtheorem{cor}[thm]{Corollary}
\newtheorem*{cor*}{Corollary}

\theoremstyle{definition}\newtheorem{defi}[thm]{Definition}
\theoremstyle{definition}\newtheorem{rem}[thm]{Remark}
\theoremstyle{definition}
\theoremstyle{definition}

\newcommand\inter{\overset{\ _{\scriptstyle\circ}}}
\numberwithin{equation}{thm}

\title{Root systems, affine subspaces, and projections
\footnote{\normalfont\copyright2021. 
{\it T\lowercase{his manuscript version is made available under the}}
CC-BY-NC-ND 4.0 
\lowercase{{\it license}\\ \url{https://creativecommons.org/licenses/by-nc-nd/4.0/}}
}}

\author{Paola Cellini}
\address{Dipartimento di Ingegneria e Geologia, Universit\`a di Chieti --
Pescara, Viale Pindaro 42, 65127 Pescara, Italy}
\email{pcellini@unich.it}

\author{Mario Marietti}
\address{Dipartimento di Ingegneria Industriale e Scienze Matematiche, Universit\`a Politecnica delle Marche, Via Brecce Bianche, 60131 Ancona,  Italy}
\email{m.marietti@staff.univpm.it}

\begin{document}

\begin{abstract}
We tackle several problems related to a finite irreducible crystallographic root system $\Phi$ in the real vector space $\mathbb E$. In particular, we study the combinatorial structure of the subsets of $\Phi$ cut by affine subspaces of $\mathbb E$ and their projections. 
As byproducts, we obtain easy algebraic combinatorial proofs of refinements of Oshima's Lemma and of a result by Kostant, a partial result towards the resolution of a problem by Hopkins and Postnikov, and new enumerative results on root systems.
\end{abstract}
\maketitle


\section{Introduction}
Let $\Phi$ be a finite irreducible crystallographic root system and let $\Pi$ be a set of simple roots of $\Phi$.
Given a subset $S$ of $\Pi$ and $\beta \in \Phi$, with $\beta=\sum\limits_{\alpha\in \Pi}c_\alpha(\beta)\alpha$,  we set
\begin{eqnarray*}
\Phi_{S,\beta}&=&\Big\{\gamma= \sum\limits_{\alpha\in \Pi}c_\alpha(\gamma) \alpha\in \Phi : c_\alpha(\gamma)=c_\alpha(\beta)\ \text{ for all } \alpha\in S\Big\},\\
\Phi_{S,\mathbb Z \beta}&=&\Big\{\gamma= \sum\limits_{\alpha\in \Pi}c_\alpha(\gamma) \alpha\in \Phi: \exists k \in \mathbb Z \text{ with }c_\alpha(\gamma)=k  c_\alpha(\beta)\ \text{ for all }\alpha\in S\Big\}.\\
\end{eqnarray*}

If $\supp(\be)\cap S=\emptyset$, the set $\Phi_{S,\beta}$ coincides with $\Phi_{S,\mathbb Z \beta}$ and is the standard parabolic subsystem generated by $\Pi\smeno S$.
So the nontrivial case happens when $\supp(\be)\cap S\neq\emptyset$.

By definition, $\Phi_{S,\be}$ is the intersection of $\Phi$ with the affine subspace  that contains $\be$ and has $\gen(\Pi\smeno S)$ as directional linear subspace. 
The subset $\Phi_{S,\mathbb Z \beta}$ is, in fact, a root subsystem of $\Phi$ and is a union of strata each obtained by intersecting  $\Phi$ with an affine subspace parallel to $\gen(\Pi\smeno S)$. 

In this paper, we study several problems concerning the subsets $\Phi_{S,\be}$ and $\Phi_{S,\mathbb Z \beta}$  and their applications. 

Section~\ref{sec:notation} provides the notational conventions used throughout the paper.

In Section~\ref{sec:rootpolytopes}, we give a uniform case-free description of the root polytope of an arbitrary (possibly reducible) root system. The {\em root polytope} of a root system is the convex hull of all roots and, in the irreducible case, it is studied in \cite{CM2}, \cite{CM1}, \cite{CM3}. The results in this section are needed in Section~\ref{sec:h-p}.

In Section~\ref{sec:Squalunque}, we study the root subsystem $\Phi_{S,\mathbb Z \beta}$ and find two of its simple systems that provide in a natural way the existence of the minimum and the maximum of $\Phi_{S, \be}$ viewed as a subposet of the root poset. 

As consequences of the results in this section,  we obtain easy case-free algebraic combinatorial proofs of a result by Kostant,  of a result by Oshima, and some refinements. Kostant's result, which is presented in \cite[Section~2]{J}, asserts the irreducibility of a certain representation of the semisimple subalgebra $\mathfrak g_{\Pi\smeno S}$  associated with $\Pi\smeno S$ (see Remark~\ref{rem:kostant}).  Oshima's result \cite[Lemma~4.3]{O}, which is Corollary~\ref{cor:transitivo} in this paper, has received attention recently for its applications and is referred to  as \lq\lq Oshima's Lemma\rq\rq. The purpose of Dyer and Lehrer in \cite{D-L} is to provide an elementary self-contained proof of (a generalised version of) Oshima's Lemma. Our proof here is very elementary, self-contained, and independent of the representation theory of semisimple complex Lie algebras (see \cite{D-L}). Oshima's Lemma states that the action of $W\la\Pi\smeno S\ra$ partitions the roots in $\Phi_{S, \be}$ according to their length, hence in at most two orbits.
In the representation theory interpretation, this means that the weights of the irreducible representation studied by Kostant are divided in at most two orbits.
Indeed, the fact that the roots of maximal length constitute a single orbit follows directly from the irreducibility of this representation.  
We provide general necessary and sufficient conditions for establishing whether $\Phi_{S,\be}$ contains roots of one or two lengths (Propostion \ref{pro:min-lunga}). 
Moreover, we provide an easy combinatorial criterion on the pair $(S,\be)$ for $\Phi_{S, \be}$ to contain long roots (Proposition~\ref{pro:lacing}). 
Though this last criterion is uniform, our proof requires some case by case check.

In Sections~\ref{sec:codim1} and \ref{sec:inclusioni}, we study in more detail the case of codimension 1,  when $S$ consists of a single simple root $\alpha$. 
In this case, we obtain a partition of $\Phi$ into subsets $\Phi_{\al,k}=\Phi \cap  \{x\in \mathbb E :(  x , \widecheck\omega_\al) = k\}$ that are all intervals w.r.t. the usual partial order, except for the 0-level $\Phi_{\al,0}$, which is a possibly reducible root subsystem.
In particular, we analyse the projections of the subsets $\Phi_{\al,k}$ 
 on the hyperplane $\gen (\Pi \smeno \{\alpha \})$. 
As a consequence of that analysis, we obtain a partial result towards the resolution of a problem by Hopkins and Postnikov. In  \cite{H-P}, Hopkins and Postnikov show that the projection of the root polytope of $\Phi$ fall inside the polytope formed by expanding by a factor of $2$ the root polytope of the root subsystem generated by $ \Pi \smeno \{\alpha\}$ (see Lemma~\ref{lem:hp}
). Their proof uses the classification of root systems and finding a case-free proof is an open problem in \cite{Ho}. In Section~\ref{sec:h-p}, we give a uniform case-free argument that proves that the result of Hopkins and Postnikov is equivalent  to the fact that a certain number  $r_{\alpha}$, defined by a uniform formula (\ref{ralpha}), is smaller than $2$, for each simple root $\alpha$.  Unfortunately, this final check is done case-by-case.

In Section~\ref{sec:numerologia}, we give some enumerative results on irreducible root systems.

\par

Except for the last part of the proof of Proposition~\ref{pro:lacing} (which is not needed in the rest of the paper),  its Lemma~\ref{lem:claimclaim}, and the final check in Section~\ref{sec:h-p} discussed above, all uniform statements are proved independently of the classification of root systems.

\section{Notational conventions}
\label{sec:notation}

Let $\Phi$ be a finite (reduced, possibly reducible) crystallographic root system in the
real  vector space ${\mathbb E} =\gen_\real \Phi$ endowed with the positive definite
bilinear form $(\,\cdot\,,\cdot\,)$. We identify ${\mathbb E}$ with its dual space via its bilinear form. 
\par

We fix our  notation on the root system and its Weyl group in the following list:
\smallskip 

\renewcommand{\arraystretch}{1.2}
$
\begin{array}{@{\hskip-1.3pt}l@{\qquad}l}
\Phi^\vee&  \textrm{the set of coroots of $\Phi$}, 
\\
\alpha^\vee &  \textrm{the coroot associated with the root $\alpha \in \Phi$}, 
\\
\Pi
&  \textrm{the set of simple roots of $\Phi$}, 
\\
\Pi^\vee
&  \textrm{the set of simple coroots of $\Phi$}, 
\\
\Omega=\{\omega_\al : \al\in \Pi\}
&  \textrm{the set of fundamental weights (the dual basis of $\Pi^\vee$)},
\\
\widecheck\Omega=\{\widecheck\omega_\al : \al\in \Pi\}
&  \textrm{the set of fundamental co-weights (the dual basis of $\Pi$)},
\\
\Phi^+  &  \textrm{the set of positive roots w.r.t. $\Pi$},
\\
\mathcal{C}  &  \textrm{the fundamental chamber associated with $\Pi$}, 
\\
c_{\alpha}(\beta)  &  \textrm{the $\alpha$-coordinate of $\beta$ w.r.t. $\Pi$: $\beta=\sum_{\alpha\in  \Pi}c_{\alpha}(\beta) \alpha$},
\\
\supp(\alpha) & =\{\alpha_i \in \Pi : c_i (\alpha) \neq 0  \}, \text{ the support of $\alpha$},
\\
W   &  \textrm{the Weyl group of $\Phi$},
\\
s_{\alpha}   & \textrm{the reflection with respect to $\alpha$},
\\
W\la S\ra & \textrm{the subgroup of $W$ generated by 
$\{s_\alpha : \alpha\in S\}$ (for $S\subseteq \Phi$)},
\\
\Phi\la S\ra & \textrm{the root subsystem generated by 
$S$ (for $S\subseteq \Phi$)},
\\
\mathcal P_{\Phi}  & \textrm{the convex hull  of all roots in $\Phi$, called the {\em root polytope} of $\Phi$}.  
\end{array}$

\bigskip
For any $x\in \mathbb E$, $x$ {\it  dominant} means $x\in \mathcal C$. 

The \emph{root poset of $\Phi$} (w.r.t. the basis $\Pi$) is $\Phi^+$ with the standard partial order, i.e., for all $\alpha,\beta \in \Phi^+$ , $\alpha \leq \beta$ if and only if $\beta - \alpha$ is a nonnegative
linear combination of roots in $\Phi^+$. 
The root poset is the transitive closure of the relation  $\alpha \lhd \beta$ if and only if 
$\beta - \alpha$ is a simple root.
The root poset hence is ranked by the height function. 

In the irreducible case, we denote by $\theta$ the highest root in $\Phi$, which is the maximum of the root poset, and by $m_{\alpha}$ the coefficient  $c_{\alpha}(\theta)$, for $\alpha \in \Pi$. Furthermore,  we denote by $\widehat \Phi$ the affine root system associated to $\Phi$ and realized as follows. 
We extend $\mathbb E$ to a vector space $\mathbb E\oplus \real\delta$, and extend the bilinear form  of $\mathbb E$ to a semidefinite positive form on $\mathbb E\oplus \real\delta$ by the condition that $\real \delta$ is the kernel of the new form. 
Then we set: 

\renewcommand{\arraystretch}{1.2}
$
\begin{array}{@{\hskip-1.3pt}l@{\qquad}l}
\alpha_0 &  =-\theta+\delta,
\\ 
\widehat{\Pi} & = \{\alpha_0\}\cup\Pi, 
\\
\widehat{\Phi}  &  =\{\beta + k \delta : \beta \in \Phi, k \in \mathbb Z \} \quad \textrm{the affine root system with simple system $\widehat \Pi$},
\\
\widehat{\Phi}^+  &=\{\beta + k \delta : \beta \in \Phi^+, k \in \mathbb Z_{\geq 0} \} \cup \{-\beta + k \delta : \beta \in \Phi^+, k \in \mathbb Z_{> 0} \} \quad \textrm{the set of}\\  
{} & \quad \textrm{ (real) positive roots of $\widehat{\Phi}$
w.r.t.  $\widehat{\Pi}$}.
\end{array}$
\bigskip

For $k\in \mathbb N$, we let  $[1, k]=\{1, \ldots, k\}$.

 \section{Root polytopes of possibly reduced root systems}\label{sec:rootpolytopes} 

Recall that we denote by $\mathcal P_{\Phi}$ the convex hull  of all  roots in $\Phi$, and we call it the {\em root polytope} of $\Phi$.  
For all irreducible $\Phi$,  a uniform explicit description of the root polytope $\mathcal P_{\Phi}$ is given in \cite{CM1}.
In this section, we extend some of the results in \cite{CM1} to the reducible case.

\par

Let us first recall some results from \cite{CM1}. If $G$ is a graph with set of vertices $V$ and $U$ is a subset of $ V$, we denote by $G\smeno U$ the subgraph of $G$ induced by the subset of vertices $V\smeno U$.  

\begin{defi} \label{defi:estremali}
Let $\Phi$ be irreducible and $\widehat \Gamma$ be its extended Dynkin diagram. For any subset $I$ of $\Pi$, we say that $I$ is {\it $\widehat\Pi$-extremal} if $\widehat \Gamma\smeno I$ is connected. 
\par
We say that the root $\al$ is {\it $\widehat\Pi$-extremal} if $\{\al\}$ is.   
\end{defi}

\begin{rem}\label{rem:estremali}
\begin{enumerate}\item
If $\Phi$ is of type $\mathrm A_n$, all roots in $\Pi$ are $\widehat\Pi$-extremal. 
For all other types, the extremal roots are the leaves of the extendend Dynkin diagram other than the affine vertex.
\item
Let $\theta$ be the highest root of $\Phi$ (supposed irreducible). A subset $I$ is $\widehat\Pi$-extremal if and only if $\{-\theta\}\cup \Pi\smeno I$ is the simple system of an irreducible root subsystem of~$\Phi$.
\end{enumerate}
\end{rem}

Let $\Phi$ be irreducible with highest root $\theta=\sum_{\al\in \Pi} m_\al\al$, and set $o_\al=\widecheck\omega_\al / m_\al $ for all $\alpha\in \Pi$.

Given $I\subseteq \Pi$, we set 
\begin{equation*}
F_I=\conv\{\ga\in \Phi : (\ga,  o_\al)=1 \text{ for all }\al\in I \}
\end{equation*}
Proposition~3.3, Proposition~3.5, and Theorem~5.5 of \cite{CM1} imply the following result.
  \begin{thm}\label{thm:faccestandard}
Let $\Phi$ be irreducible.
\begin{enumerate}
\item For each $\widehat\Pi$-extremal  subset $I$, the set $F_I$ is a face of $\mathcal P_\Phi$ of dimension $\mathrm{rk}(\Phi)-
|I|$. 
\item For each subset $I$ of $\Pi$, there exists a $\widehat \Pi$-extremal subset $\ov I$ such that  $F_I=F_{\ov I}$. 
\item
For each face $F$ of $\mathcal P_\Phi$, there exist some $w\in W$ and a unique $\widehat\Pi$-extremal subset $I$ such that $F=w(F_I)$.
\end{enumerate} 
\end{thm}

Following \cite{CM1}, we call the faces of type $F_I$ the {\it standard parabolic faces}. Note that the standard parabolic  facets are exactly the faces of type $F_{\{\al\}}$ where $\al$ is a $\widehat \Pi$-estremal root.

\begin{rem}\label{rem:catenestandard}
When $\Phi$ is irreducible, for each subset $I$ of $\Pi$, there exists a root $\al\in I$ such that $\al$ is not orthogonal to $\Pi\smeno I$. 
If $I$ is $\wh \Pi$-extremal, for such a root $\al$ we have that $I\smeno\{\al\}$ is $\wh \Pi$-extremal too, if nonempty. 
This implies that any standard parabolic face either is a facet, or is contained in a standard parabolic face of higher dimension.
In particular, each standard parabolic face is contained in a standard parabolic facet. 
\end{rem}

Now we suppose that  $\Phi$ be possibly reducible.

For any face $F$ of $\mathcal P_{\Phi}$, we denote by $\inter F$ the relative interior of $F$, i.e. its interior within the affine hull of $F$. The relative interior $\inter F$ of $F$ consists of all strictly positive convex linear combinations of the vertices of $F$.
If $F$ and $F'$ are different faces, then $\inter F\cap \inter {F'}=\emptyset$.

\begin{defi}
Let $F$ be a face of $\mathcal P_{\Phi}$. 
We call $F$ {\it standard} provided that either $\inter F\cap \mathcal C\neq \emptyset$, or $F= \emptyset$.
\end{defi}

The Weyl group $W$ acts on the faces of $\mathcal P_{\Phi}$. We prove that the standard faces are a set of representatives of the orbits of this action. Given  a face $F$ with vertex set $V$, we denote by $\beta_F$ the barycenter of $V$, i.e. $\beta_F= \frac{1}{|V|} \sum_{v\in V} v$.

\begin{thm}\label{pro:unicastandard}
Let $\Phi$ be a (possibly reducible) root system and let $F$ be a face of $\mathcal P_{\Phi}$.
\begin{enumerate}
\item 
\label{11}
$F$ is standard if and only if the barycenter $\be_F$ is dominant. In this case, the stabilizer $\stab_W(F)$ of $F$ is the standard parabolic subgroup of $W$ generated by $\Pi\cap \be_F^\perp$. 
\item
\label{12}
 The orbit $W\cdot F$ of $F$ under the action of $W$ contains exactly one standard face. 
 \end{enumerate}
\end{thm}

\begin{proof}
(\ref{11})
Given $\ga\in \Phi$, we denote by $s_\ga$ the reflection through $\ga$ and by $H_{\ga}$ the hyperplane fixed by $s_\ga$.

For each $w\in W$, if $w(F)\neq F$ then $\inter F\cap w(\inter F)=\emptyset$. 
In particular, for each $\ga\in \Phi$, since $H_\ga\cap \inter F$ is fixed by $s_\ga$, if $s_\ga(F)\neq F$ then $H_\ga\cap \inter F=\emptyset$, and $\inter F$ lies in one of the two open half spaces determined by $H_\ga$. 
If, in addition, $\ga>0$ and $\inter F\cap \mathcal C\neq\emptyset$, then $(x, \ga) > 0$ for all $x\in \inter F\cap \mathcal C$, and hence $(x, \ga)\geq 0$ for all $x\in F$.
\par
The stabilizer $\stab_W(F)$ of $F$ coincides with the stabilizer of the baricenter $\beta_F$,  hence $\stab_W(F)$ is the parabolic subgroup of $W$ generated by $\{s_\ga: (\ga, \be_F)=0\}$. 
We have just seen that, if $F$ is standard, then $(\ga, \be_F)> 0$ for all $\ga\not\in\stab_W(F)$, hence $\be_F\in \mathcal C$. 
\par
Conversely, if $\be_F\in \mathcal C$, then $F$ is standard by definition, since $\be_F\in \inter F$. 
\par
Furthermore, if $\be_F\in \mathcal C$, then $\stab_W(F)$ is a standard parabolic subgroup, namely, $\stab_W(F)=W\la \Pi\cap \be_F^\perp \ra$.
\par
(\ref{12})
Since $\mathcal C$ is a fundamental domain for $W$, there is one and only one face in $W\cdot F$ with dominant barycenter, hence $W\cdot F$ contains exactly one standard face by  (\ref{11}). 
\end{proof}

As a direct consequence of Theorem~\ref{pro:unicastandard}, we have the following result in the irreducible case.
\begin{cor}
\label{cor:lostesso}
Let $\Phi $ be irreducible. The standard faces of $\mathcal P_{\Phi}$ are exactly the standard parabolic faces.
\end{cor}
\begin{proof}
 Since the barycenters of the standard parabolic faces are dominant (see \cite[Proposition~4.4]{CM1}), the statement follows by Theorem~\ref{pro:unicastandard}.
\end{proof}

As we see below, each dominant point belongs to the cone of some standard face.

\begin{pro}\label{pro:semistandard}
Let $\Phi$ be a (possibly reducible) root system and let $F$ be a face of $\mathcal P_{\Phi}$ such that $F\cap \mathcal C\neq \emptyset$. 
Then there is a standard subface $F'$ of $F$ such that $F\cap \mathcal C=F'\cap \mathcal C$.
\end{pro}

\begin{proof}
If $F$ is a vertex of $\mathcal P_{\Phi}$, then $F$ is standard.
If $F$ is not standard, then there is a proper subface $F'$ of $F$ such that $F\cap \mathcal C=F'\cap \mathcal C$. 
The assertion follows by induction on the dimension of $F$.
\end{proof}

Proposition~\ref{pro:semistandard} directly implies  the following result.
 
\begin{cor}\label{cor:unionefaccestandard}
Let $\Phi$ be a (possibly reducible) root system and let $\partial \mathcal P_{\Phi}$ be the boundary of  $\mathcal P_{\Phi}$, i.e. the union of its faces. 
Then $\mathcal C\cap \partial \mathcal P_{\Phi}$ is contained in the union of the standard faces. 
In particular, for each nonzero $x\in \mathcal C$, there exist a positive real number $r$ and a standard face $F$  such that $r x\in F$. 
\end{cor}

We now determine explicitly the standard faces in the case $\Phi $ is any (possibly reducible) root system.

Let $\Phi_1, \dots, \Phi_k$ be the irreducible components of $\Phi$, $\mathbb E_i=\gen(\Phi_i)$, $\Pi_i=\Pi\cap \Phi_i$, $\mathcal C_i=\mathcal C\cap \mathbb E_i$, and let $\mathcal P_i=\mathcal P_{\Phi} \cap \mathbb E_i$ be the root polytope of $\Phi_i$. 
We say that $x\in \mathbb E$ is $\Phi_i$-dominant provided that $x\in\mathcal C_i$.

\par
Let $F$ be a face of $\mathcal P_{\Phi}$ and  $$F_i=F\cap \mathbb E_i.$$  
If $F=\mc P_{\Phi}$, then $F_i=\mc P_i$ for all $i\in [1, k]$. 
Conversely, if $F\neq \mc P_\Phi$, then  $F_i\neq \mc P_i$ for all $\in [1,k]$. 
By a {\em proper face}, we intend a possibly empty face different from the whole polytope (this definition slightly differs from the usual one but  is more convenient for our purposes). 

\begin{pro}\label{lem:faccestandard}
Let $\Phi$ be a (possibly reducible) root system. 
\begin{enumerate}
\item
\label{item: un} The map sending a face $F$ to $(F_1, \ldots, F_k)$, where $F_i=F\cap \mathbb E_i$,  is a bijection between the proper faces of $\mc P_{\Phi}$ and the $k$-tuples $\ov F_1, \dots, \ov F_k$, where $\ov F_i$ is  a proper face of  $\mc P_i$ for every $i\in[1, k]$.
\item 
\label{item: du}
A proper face $F$ has dimension 
$$\dim F=k-1+\sum\limits_{i=1}^k \dim F_i.$$
\item 
\label{item: tr}A face $F$ of  $\mathcal P_{\Phi}$ is standard if and only if $F_i$ is standard for all $i\in [1,k]$.
\end{enumerate}
\end{pro}

\begin{proof} (\ref{item: un})
Let  $F$ be a proper face of $\mc P_{\Phi}$, and denote by $H_F$ a supporting hyperplane for $F$, i.e., a hyperplane  such that $F= \mc P_{\Phi} \cap H_F$ and $\mc P_{\Phi}$ is contained in one of the two halfspaces cut by $H_F$. 
Since $\Phi$ is symmetric and spans $\mathbb E$, the hyperplane $H_F$ does not contain the origin, hence $H_F\cap \mathbb E_i$ is either empty, or a hyperplane in $\mathbb E_i$.  
Since $F_i=F\cap \mathbb E_i=\mathcal P_{\Phi}\cap H_F\cap \mathbb E_i=\mathcal P_i\cap(H_F\cap \mathbb E_i)$, it follows that $F_i$ is a proper (possibly empty) face of~$\mathcal P_i$. 
Since $F=\conv(F_1\cup\dots\cup F_k)$, the map in the statement is injective.
Conversely, for any choice of proper faces $\ov F_1, \dots, \ov F_k$ of $\mc P_1, \dots, \mc P_k$, respectively, consider  the convex hull $F$ of $\ov F_1\cup\dots\cup \ov F_k$. Clearly, $F \cap \mathbb E_i=\ov F_i$, for every $i\in[1, k]$.  We need to show that  $F$ is a proper face of $\mc P_{\Phi}$. 
Indeed, since the supporting hyperplanes of the $\ov F_i$ do not contain the origin, there exist $y_i\in \mathbb E_i$ such that  $(y_i, x)\leq 1$  for all $x\in \mc P_i$, and $\ov F_i=\{x\in \mc P_i: (x, y_i)=1\}$, for every $i\in[1, k]$. 
If we define $y=y_1+\dots+y_k$, we obtain $(y, x)\leq 1$ for all $x\in \mc P_{\Phi}$ and $(y, x)=1$ if and only if $x\in F$.

\par 
(\ref{item: du})
Since a face of a polytope is the convex hull of the vertices of the polytope that belong to the face,   
$F$ is the convex hull of $F\cap \Phi$, and $F_i$ is the convex hull of $F\cap \Phi_i$. 
If $F$ is proper and nonempty, then $\dim F=\dim(\gen(F\cap \Phi))-1$ and $\dim F_i=\dim (\gen (F\cap \Phi_i))-1$
The dimension formula follows by the preceding formulas (which hold also for the empty face, since it has dimension $-1$ by definition). 

\par
(\ref{item: tr})
The set of vertices of $F$ is the union of the sets of vertices of the subfaces $F_i$.
Hence, the barycenter $\be_F$ of $F$ is a positive linear combination of the barycenters $\be_{F_i}$ of the nonempty faces $F_i$, and $\be_F$ is dominant if and only if all these  barycenters $\be_{F_i}$ are $\Phi_i$-dominant. 
This implies the claim by Theorem~\ref{pro:unicastandard}.
\end{proof}

Now let $F$ be standard. For $i\in [1,k]$,   denote by $\theta_i$ the highest root of $\Phi_i$. For $\alpha \in \Pi_i$, define $m_{\alpha}$ by the condition $\theta_i=\sum_{\al\in \Pi_i} m_\al\al$, and set  $o_\al=\widecheck\omega_\al / m_\al$.
By Theorem~\ref{thm:faccestandard}, Corollary~\ref{cor:lostesso}, and Proposition~\ref{lem:faccestandard}, for every nonempty $F_i=F\cap \mathbb E_i$, there exists a unique $\widehat \Pi_i$-extremal subset $\Pi_i(F)$ of $\Pi_i$  such that:
$$F_i=\{x\in \mathcal P_i: (x, o_\al)=1\text{ for all } \al\in \Pi_i(F)\}.$$
We let 
\begin{equation*}\label{eq:defoF} 
o_i(F)=\frac{1}{|\Pi_i(F)|}\sum\limits_{\al\in \Pi_i(F)} o_\al.
\end{equation*}  
For $F_i=\emptyset$, we set $o_i(F)=0$.
Since $(x, o_\al)\leq 1$ for all $\al\in \Pi_i$ and $x\in \mc P_i$, 
we have
\begin{equation}\label{eq:oF}
(x, o_i(F))\leq 1, \text{ for all } x\in \mc P_i, \text{ and } F_i=\{x\in \mathcal P_i: (x, o_i(F))=1\}.
\end{equation}
The following result holds.
\begin{pro}
Let $\Phi$ be a (possibly reducible) root system. If $F$ is a standard face of  $\mathcal P_{\Phi}$, then $\left(x, \sum\limits_{i=1}^k o_i(F)\right)\leq 1$ for all $x\in \mc P_{\Phi}$ and  
$$
F=\left\{x\in \mathcal P_{\Phi} : \left(x, \sum\limits_{i=1}^k o_i(F)\right)=1 \right\}.
$$
\end{pro}

\begin{proof}
Let $x\in \mc P_{\Phi} $. From any expression of $x$ as a convex linear combination of the roots in $\Phi$, we easily obtain a subset $J$ of $[1, k]$ and nonzero elements $x_i\in \mc P_i$,  for all $i\in J$, such that $x$ is a  strictly positive convex linear combination of such elements $x_i$.
By \eqref{eq:oF} above,   $(x_i, o_i(F))\leq 1$ for all $i\in J$, hence $\left(x, \sum\limits_{i=1}^k o_i(F)\right)=\left(x, \sum\limits_{i\in J} o_i(F)\right)\leq 1$. 
Moreover, $\left(x, \sum\limits_{i\in J} o_i(F)\right)=1$ if and only if $(x_i, o_i(F))=1$ for all $i\in J$, hence if and only if $x_i \in F_i$ for all $i\in J$: this  is equivalent to $x\in F$.
\end{proof}

In the following corollary, we isolate the special case of facets. 
The proof is a direct application of the preceding results. 
\begin{cor}\label{cor:facetsstandard}
Let $\Phi$ be a (possibly reducible) root system and let $F$ be a facet of $\mathcal P_{\Phi}$. Then $F_i\neq \emptyset$ and $\dim F_i=\dim \mathbb E_i-1$,  for all $i\in[ 1,k]$.
\par
In particular, $F$ is standard if and only $F_i$ is a standard facet of $\mathcal P_i$. 
In this case, for all $i\in [1,k]$, the set $\Pi_i(F)$ is a singleton containing a $\widehat\Pi_i$-extremal root $\al_i\in \Pi_i$, hence $o_i(F)=o_{\al_i}$ and the supporting hyperplane of $F$ is 
$$\{x\in \mathbb E: (x, o_{\al_1}+\dots+ o_{\al_k})=1\}.$$
Conversely, for each $k$-tuple $(\al_1, \dots, \al_k)\in \Pi_1\times\dots\times \Pi_k$ such that $\al_i$ is  $\widehat\Pi_i$-extremal  for all $i\in [1, k]$, the hyperplane
$$\{x\in \mathbb E: (x, o_{\al_1}+\dots+ o_{\al_k})=1\}$$ 
supports a standard facet of $\mathcal P_{\Phi}$.
\end{cor}

Corollary~\ref{cor:facetsstandard} and Remark~\ref{rem:catenestandard} together imply that each standard face is a subface of a standard facet. Hence, by Corollary~\ref{cor:unionefaccestandard}, we obtain the following result. 

\begin{cor}\label{cor:inclusionedominanti}
Let $x$ be a dominant element in $\mathbb E$. 
Then $\conv(W\cdot x)\subseteq \mathcal P_{\Phi}$ if and only if
$$(x, o_{\al_1}+\dots+ o_{\al_k})\leq 1,$$ 
for each $k$-tuple $(\al_1, \dots, \al_k)\in \Pi_1\times\dots\times \Pi_k$ such that $\al_i$ is  $\widehat\Pi_i$-extremal  for all $i\in [1, k]$.
\end{cor}

\section{Root systems and affine subspaces}\label{sec:Squalunque}

In this section, we study the intersections  of an irreducible root system  with certain  affine subspaces. 

Let $\Phi$ be irreducible.

Given a subset $S$ of $\Pi$ and $\beta \in \Phi$, with $\beta=\sum\limits_{\alpha\in \Pi}c_\alpha(\beta)\alpha$,  we set
\begin{eqnarray*}
\Phi_{S,\beta}&=&\Big\{\gamma= \sum\limits_{\alpha\in \Pi}c_\alpha(\gamma) \alpha\in \Phi : c_\alpha(\gamma)=c_\alpha(\beta)\ \text{ for all } \alpha\in S\Big\},\\
\Phi_{S,\mathbb Z \beta}&=&\Big\{\gamma= \sum\limits_{\alpha\in \Pi}c_\alpha(\gamma) \alpha\in \Phi: \exists k \in \mathbb Z \text{ with }c_\alpha(\gamma)=k  c_\alpha(\beta)\ \text{ for all }\alpha\in S\Big\}.\\
\end{eqnarray*}

\begin{rem}
\label{senza loss}
Let $\beta \in \Phi$. If $\Phi$  were reducible and $\Psi$ were the irreducible component of $\Phi$ containing $\beta$, then $\Phi_{S,\beta}=\Psi_{S,\beta}$ and $\Phi_{S,\mathbb Z  \beta}=\Psi_{S,\mathbb Z \beta}$ for all nonempty $S$. Thus, there is no loss of generality in our assumption   that $\Phi$  be irrreducible.
\end{rem}

Recall that $\delta$ denotes the indivisible imaginary root of the affine root system $\widehat \Phi$ (see Section~\ref{sec:notation}).
\begin{pro}\label{pro:sottosistema}
Let $S\subseteq \Pi$ and $\beta\in \Phi^+$ with $\supp(\beta)\cap S\neq \emptyset$.

\begin{enumerate}
\item
\label{item: uno}
$\Phi_{S,\mathbb Z \beta}$ is a root subsystem of $\Phi$.
\item
\label{item: due}
 $\Phi_{S,\beta}$ is an interval of the root poset, i.e. it has minimum and maximum and contains every root in between. Moreover,  both $$\{\min\Phi_{S,\beta}\}\cup(\Pi\smeno S) \quad \text{and} \quad
\{-\max\Phi_{S,\beta}\}\cup(\Pi\smeno S)$$ 
are simple systems of $\Phi_{S,\mathbb Z \beta}$;
in particular, $\min\Phi_{S,\beta}$ and $\max\Phi_{S,\beta}$ have the same length.
\item 
\label{item: tre}
Let $\Psi=\Phi_{S,\mathbb Z \beta}$ and choose  $\{\gamma\} \cup (\Pi \smeno S)$ as a set of simple roots of $\Psi$, with either $\gamma= \min\Phi_{S,\beta}$ or $\gamma= - \max\Phi_{S,\beta}$. Then either $$\Phi_{S,\beta}=\Psi_{\{\gamma\},\gamma} \quad \textrm{ or } \quad  \Phi_{S,\beta}= - \Psi_{\{\gamma\},\gamma},$$
respectively.
\end{enumerate}
\end{pro}
\begin{proof}
(\ref{item: uno}) The assertion is clear.
\par
(\ref{item: due})
Let $U=\gen\left( (\Pi\smeno S)\cup\{\beta+\delta\}\right)$ and 
$\widehat \Phi_U=\widehat\Phi\cap U$.  
Then $\Phi_{S,\beta}+\delta\subseteq \widehat\Phi_U$. \par
Since $\supp(\beta)\cap S$ is nonempty, $\real \delta\cap U=\{0\}$, hence the restriction  of $(\cdot\, , \cdot)$ to $U$ is positive definite and $\widehat \Phi_U$ is a finite dimensional crystallographic root system, naturally isomorphic to 
$\Phi_{S,\ganz\beta}$ through the natural projection $E+\real \delta\to E$.
We can choose \hbox{$\widehat\Phi_U\cap \widehat \Phi^+$} as a positive system for $\widehat\Phi_U$. 
The simple system determined by this choice is  \hbox{$(\Pi\smeno S)\cup \{\gamma+\delta\}$}, for a certain $\gamma\in \Phi_{S,\beta}$ with $\gamma\leq \beta$. 
Indeed, since all roots in  $\Phi_{S,\beta}+\delta$ must be non-negative linear combination of the basis,  $\gamma=\min \Phi_{S,\beta}$. 
\par
If we make the above construction with ${-\beta}+\delta$ in place of $\beta+\delta$, 
we obtain, as a simple system for $\widehat\Phi_U$,  $(\Pi\smeno S)\cup \{-\gamma+\delta\}$, with $\gamma=\max \Phi_{S,\beta}$. 

By definition of $\Phi_{S,\beta}$, every root  between $\min \Phi_{S,\beta}$ and $\max \Phi_{S,\beta}$ in the root poset belongs to $\Phi_{S,\beta}$.
\par
Finally, it is clear that, if  $(\Pi\smeno S)\cup \{\gamma+\delta\}$ is a basis of $\widehat\Phi_U$, then $(\Pi\smeno S)\cup \{\gamma\}$ is a basis of 
$\Phi_{S,\mathbb Z \beta}$. 
\par
(\ref{item: tre})
The assertion follows by (\ref{item: uno}) and (\ref{item: due}). 
\end{proof}

\begin{rem}
\label{rem:kostant}
Proposition~\ref{pro:sottosistema} implies a representation theoretical result due to Kostant, which is presented in \cite[Section~2]{J}.
We briefly recall this result.
Let $\mathfrak g$ be a simple Lie algebra with root system $\Phi$, and $\mathfrak p_{\Pi\smeno S}$ be the standard parabolic subalgebra of $\mathfrak g$ corresponding to $\Pi\smeno S$.
Then
$$\mathfrak p_{\Pi\smeno S}=\mathfrak g_{\Pi\smeno S}\oplus\mathfrak h_{(\Pi\smeno S)^\perp}\oplus \mathfrak n^+_{\Phi\smeno \Phi(\Pi\smeno S)},$$ where $\mathfrak g_{\Pi\smeno S}$ is the semisimple subalgebra generated by the root spaces $\mathfrak g_\al$ with $\al\in \Pi\smeno S$, $\mathfrak h_{(\Pi\smeno S)^\perp}$ is defined as $\{h \in \mathfrak h : \alpha(h)=0, \,  \forall  \alpha \in \Pi\smeno S \}$ and is the center of the reductive subalgebra $\mathfrak g_{\Pi\smeno S}\oplus\mathfrak h_{(\Pi\smeno S)^\perp}$, and  $\mathfrak n^+_{\Phi\smeno \Phi(\Pi\smeno S)}=\bigoplus\limits_{\ga \in \Phi^+\smeno \Phi(\Pi\smeno S)} \mathfrak g_{\ga}$.
The adjoint representation of $\mathfrak g$ restricts to a representation of $\mathfrak g_{\Pi\smeno S}\oplus\mathfrak h_{(\Pi\smeno S)^\perp}$ on the nilpotent summand $\mathfrak n^+_{\Phi\smeno \Phi(\Pi\smeno S)}$.
For any $\beta\in \Phi^+$, if $\beta=\sum\limits_{\alpha\in \Pi}c_\alpha\alpha$, the action of  $\mathfrak h_{(\Pi\smeno S)^\perp}$ on  $\mathfrak g_\beta$ depends only on $\sum\limits_{\alpha\in S}c_\alpha\alpha$, hence,  $\bigoplus\limits_{\ga \in \Phi_{S, \be}} \mathfrak g_{\ga}$, which we denote  by $\mathfrak n_{\Phi_{S, \beta}}$,  is a $\mathfrak h_{(\Pi\smeno S)^\perp}$-isotypic component.  
Moreover, $\mathfrak n_{\Phi_{S, \beta}}$ is stable under the action of $\mathfrak g_{\Pi\smeno S}$.
Kostant's results is that,  whenever $\supp(\beta)\cap S\neq  
\emptyset$, the subrepresentation of $\mathfrak g_{\Pi\smeno S}$  on   
$\mathfrak n_{\Phi_{S, \beta}}$ is irreducible.
The whole set of weights of this irreducible representation is $\Phi_{S, \be}$, hence its unique highest weight is the maximum of $\Phi_{S, \be}$. We show that the set of weights $\Phi_{S, \be}$ is partioned into one or two $W\la\Pi \smeno S\ra$-orbits according to Propositions~\ref{pro:come-beta},~\ref{pro:min-lunga}, and~\ref{pro:lacing} (see also Remark~\ref{rem:irriducibile}).

We also note that, by a general result of representation theory (see \cite[p. 204] {F-H}), the set of weights  ${\Phi_{S, \beta}}$  is contained in the convex hull of the orbit of the highest weight under the action of  the Weyl group of $\mathfrak g_{\Pi\smeno S}$, hence $$ \conv(\Phi_{S, \be}) = \conv \big( W\la\Pi \smeno S\ra\cdot \max\Phi_{S, \be} \big).
$$
\end{rem}
\begin{rem}

Let $\alpha \in \Pi$ and $S=\{\alpha\}$. Then $\Phi_{S,\mathbb Z \alpha}=\Phi$ and $\{\min\Phi_{S,\alpha}\}\cup(\Pi\smeno S)= \Pi$. By Proposition~\ref{pro:sottosistema}, also $\{-\max\Phi_{S,\alpha}\}\cup(\Pi\smeno S)$ is a  simple system of $\Phi$.
When $m_{\alpha}=1$, $\max\Phi_{S,\alpha}=\theta$ and, since the affine simple root $\alpha_0$ satisfies $(\alpha_0,x)=(-\theta,x)$ for all $x\in \mathbb E$, we reobtain  the following well-known fact:

{\em the coefficient  $m_{\alpha}$ is equal to  $1$ if and only if the Dynkin diagram of $\Pi$ is equal to the subdiagram  $\Pi \cup \{\alpha_0\} \smeno \{\alpha\} $ of the extended Dynkin diagram of $\Phi$.}
\end{rem}

\bigskip 
Recall that $\widecheck\omega_\al$, for each $\al\in \Pi$, denotes the fundamental co-weight corresponding to $\al$, i.e. the element of $\mathbb E$ defined by the conditions $(\widecheck\omega_\al, \al)=1$ and $(\widecheck\omega_\al, \al')=0$ for $\al'\in \Pi\smeno \{\al\}$. 
Let 
$$\mathbb A=\gen(\Pi\smeno S), \qquad\text{ so that } \qquad \mathbb A^\perp=\gen(\widecheck\omega_\al : \al\in S),$$ 
and denote by $\pi_\mathbb  A$ and $\pi_{\mathbb A^\perp}$ the orthogonal projections from $\mathbb E$ onto $\mathbb  A$ and $\mathbb  A^\perp$, respectively.

\begin{lem}\label{lem:proiezione}
For all $\ga \in \Phi_{S, \be}$ and $w\in W\la\Pi\smeno S\ra$, 
\begin{enumerate}
\item\label{item:trasla}
$\pi_{\mathbb A^\perp}(\ga)=\pi_{\mathbb A^\perp}(\be)$, hence $\pi_\mathbb  A(\ga)=\ga - \pi_{\mathbb  A^\perp}(\be)$, 
\item\label{item:commuta} 
$\pi_\mathbb  A \circ w(\ga)= w\circ \pi_\mathbb  A (\ga)$.
\end{enumerate}
\end{lem}

\begin{proof}
The equalities in (1) follow directly from the definitions.  
Equality (2) follows easily from (1), since $W\la\Pi\smeno S\ra$ stabilizes $\mathbb  A$ and fixes $\mathbb  A^\perp$ pointwise. 
\end{proof}

Given a subset $T$ of $\Pi$, we say that $v\in \mathbb E$ is {\it $T$-dominant} provided that $(v,\alpha)\geq 0$ for all $\alpha\in T$. 
We note that, for all $v\in \mathbb E$, $v$ is $(\Pi\smeno S)$-dominant if and only if $\pi_\mathbb A(v)$ belongs to the fundamental chamber of the possibly reducible root system $\Phi\la\Pi\smeno S\ra$.

\begin{pro}\label{pro:intervallo-gen}
Let  $S\subseteq \Pi$ and $\beta\in \Phi^+$, with $\supp(\beta)\cap S\neq \emptyset$. Then
\begin{enumerate} 
\item
\label{item:1-intervallo-gen} $\max \Phi_{S,\beta}$ and $- \min \Phi_{S,\beta}$ are both  $(\Pi\smeno S)$-dominant; 
\item
\label{item:3-intervallo-gen}
 if  $\Phi_{S,\beta}$ contains some long roots, then  $\max \Phi_{S,\beta}$ is the unique long $(\Pi\smeno S)$-dominant root; otherwise, $\max \Phi_{S,\beta}$ is the unique (short) $(\Pi\smeno S)$-dominant root.
\end{enumerate}
\end{pro}

\begin{proof}
Assertion (\ref{item:1-intervallo-gen}) is implicit in Proposition \ref{pro:sottosistema} (\ref{item: due}), and may also be seen directly 
since, for each $\ga\in \Phi$, if $(\gamma, \alpha)<0$ for some $\alpha\in \Pi\smeno S$, then $\gamma+\al\in \Phi$. 
\par

Let us prove (\ref{item:3-intervallo-gen}). We first show that, for each $\ga\in \Phi_{S,\beta}$, there exists  a $(\Pi\smeno S)$-dominant root of the same length as $\ga$.
The parabolic subgroup $W\la\Pi\smeno S\ra$ acts on the set of long (short) roots in $\Phi_{S,\beta}$. By Lemma~\ref{lem:proiezione},  the orthogonal projection $\pi_\mathbb A:\mathbb E\to \mathbb A$ restricts to a bijection from the $W\la\Pi\smeno S\ra$-orbit of $\ga$ to the  $W\la\Pi\smeno S\ra$-orbit of $\pi_\mathbb A(\ga)$. 
Since $W\la\Pi\smeno S\ra\cdot \pi_\mathbb A(\ga)$ contains a unique element belonging to the fundamental chamber of $\Phi(\Pi\smeno S)$, the orbit $W\la\Pi\smeno S\ra\cdot \ga$ contains a unique $(\Pi\smeno S)$-dominant root.

Next, we prove that, if $\gamma$ is a $(\Pi\smeno S)$-dominant root in $\Phi_{S, \beta}$ other than $\max \Phi_{S, \beta}$, then   $\max \Phi_{S, \beta}$ is a long root and  $\gamma$ is short. Since $\gamma<\max \Phi_{S, \beta}$, we have  $\max \Phi_{S, \beta}- \gamma \in \mathbb N ( \Pi \smeno S)$, and hence  $(\max \Phi_{S, \beta},\gamma^\vee)\geq(\gamma, \gamma^\vee)=2$, since $(\alpha, \gamma^\vee) \geq 0$ for all $\alpha \in \Pi \smeno S$. This implies that  $\max \Phi_{S, \beta}$ is long and $\gamma$ is short.  

Hence, if $\Phi_{S, \beta}$ contains either only long roots or both long and short roots, then $\max \Phi_{S, \beta}$ is the unique long $(\Pi\smeno S)$-dominant root. Moreover, if $\max \Phi_{S, \beta}$ is short, then it is the unique short $(\Pi\smeno S)$-dominant root in $\Phi_{S, \beta}$, and there is no long root in $\Phi_{S, \beta}$. 
\end{proof}

As a consequence, we obtain an algebraic combinatorial short proof of Oshima's Lemma (see \cite{O}), which is case-free and independent of the representation theory of semisimple complex Lie algebras (see \cite{D-L}). 
\begin{cor}
[Oshima's Lemma]\label{cor:transitivo}
Let $S\subseteq \Pi$ and $\beta\in \Phi^+$, with $\supp(\beta)\cap S\neq \emptyset$. Then
\begin{enumerate}
\item
\label{unouno}  if $\Phi_{S,\beta}$ contains a long (short) root, then it contains exactly one $(\Pi\smeno S)$-dominant long (short) root;
\item
\label{duedue} 
 $W\la\Pi\smeno S\ra$ is transitive on the set of long (short) roots in $\Phi_{S, \beta}$.
 \end{enumerate}
\end{cor}
\begin{proof}
Let us prove (\ref{unouno}). 
By Proposition~\ref{pro:intervallo-gen}(\ref{item:3-intervallo-gen}), if $\Phi_{S, \beta}$ contains some long roots, then it contains exactly one $(\Pi\smeno S)$-dominant long root. 

It remains to show that, if $\Phi_{S, \beta}$ contains some short roots, then it contains exactly one $(\Pi\smeno S)$-dominant short root. We prove it by duality. 

For a subset $T$ of $ \Phi$, let $T^\vee=\{\al^\vee:\al\in T\}$. 
Let $\be_s$ be any short root in $\Phi_{S, \beta}$. Then $\ga$ is a short root in $\Phi_{S, \beta}$ if and only if $\ga^\vee$ is a long root in $\Phi^\vee_{S^\vee,\be_s^\vee}$, and $\ga$ is  $(\Pi\smeno S)$-dominant root if and only if $\gamma^\vee$ is $(\Pi^\vee\smeno S^\vee)$-dominant. Hence, we are done.

Part (\ref{2}) is, in fact, equivalent to Part (\ref{unouno}) since, as already noticed,  Lemma~\ref{lem:proiezione} implies that  the orthogonal projection $\pi_\mathbb A:\mathbb E\to \mathbb A$ restricts to a bijection from the $W\la\Pi\smeno S\ra$-orbit of each $\ga\in \Phi_{S, \beta}$ to the  $W\la\Pi\smeno S\ra$-orbit of $\pi_\mathbb A(\ga)$. 
\end{proof}

\begin{rem}
If $\Phi_{S, \beta}$ contains some  short roots, then the unique short $(\Pi\smeno S)$-dominant root is also the unique maximal short root in $\Phi_{S, \beta}$. Indeed, if $\gamma\in \Phi_{S, \beta}$ is not $(\Pi\smeno S)$-dominant, for $\alpha\in \Pi\smeno S$ with $(\gamma, \alpha)<0$, we have that  $s_\alpha(\gamma) >\gamma$, so $\gamma$ cannot  be maximal among the short roots.
\end{rem}

\begin{rem}\label{rem:irriducibile}
In Remark~\ref{rem:kostant}, we have seen that $\max \Phi_{S, \be}$ is the highest weight of an irreducible representation of the semisimple Lie algebra associated to $\Phi\la\Pi \smeno S\ra$.
Corollary~\ref{cor:transitivo} shows that the action of $W\la\Pi \smeno S\ra$ on the weights of this representation has one orbit for each length of the roots $\Phi_{S, \be}$, hence at most two orbits. 
In particular, in the simply laced case there is  a single orbit action. 
By Proposition~\ref{pro:intervallo-gen}, there is a single orbit action also when $\max \Phi_{S, \be}$ is short. 
\end{rem}

In the rest of the section, we study when $\Phi_{S, \be}$ contains roots of one or two lengths.

Recall that the lacing number of $\Phi$ is  the maximal number of edges connecting two vertices of the Dynkin diagram of $\Phi$, i.e. $\max\{ (\al,\al) / (\be,\be) : \alpha, \beta \in \Phi \}$. 
As usual, in simply laced root system, we consider all roots as long.
The following result is an exercise of Bourbaki's \cite{Bou}.
We include a proof for completeness.

\begin{lem}\cite[Exercise 20, Chapter VI, \S1]{Bou}\label{lem:esercizio-bou}
Let $\ga= \sum_{\alpha \in \Pi} c_{\alpha}\alpha \in \Phi$.  
Then, $\ga $ is long 
if and only if, for all short $\alpha$ in $\Pi$,  the coefficient $ c_{\alpha}$ is a multiple of the lacing number of $\Phi$.
\end{lem}

\begin{proof}
Let $\ga$ be a long root. 
The coroot $\ga^\vee$ satisfies 
$$ \ga^\vee = \frac{2 \ga}{(\ga, \ga)}= \sum_{\alpha \in \Pi_\ell} \frac{2 }{(\ga, \ga)}c_{\alpha}\alpha + \sum_{\alpha \in \Pi_{s}} \frac{2 }{(\ga, \ga)}c_{\alpha}\alpha=\sum_{\alpha \in \Pi_\ell} c_{\alpha}\alpha^\vee + \sum_{\alpha \in \Pi_{s}} \frac{(\alpha,\alpha )}{(\ga, \ga)}c_{\alpha}\alpha^\vee,$$
where $\Pi_{\ell}$ and $\Pi_s$ denote, respectively, the set of simple long and short roots.
Now,  for all short root $\alpha$, the ratio $(\alpha,\alpha ) / (\ga, \ga)$ coincides with the inverse of the lacing number of $\Phi$; since the coefficients of $\ga^\vee$ must be integers, $ c_{\alpha}$ is a multiple of the lacing number,  for each short $\alpha \in \Pi$.

\par
Conversely, let $\gamma=\sum_{\alpha \in \Pi}d_{\alpha}\alpha$ be a root such that $ d_{\alpha}$  is a multiple of the lacing number of $\Phi$, for all short $\alpha \in \Pi$. We prove that $\gamma$ is long by contradiction. Suppose that it is short. The coroot $\gamma^\vee$ satisfies 
$$ \gamma^\vee = \frac{2 \gamma}{(\gamma, \gamma)}= \sum_{\alpha \in \Pi_\ell} \frac{2 }{(\gamma, \gamma)}d_{\alpha}\alpha + \sum_{\alpha \in \Pi_{s}} \frac{2 }{(\gamma, \gamma)}d_{\alpha}\alpha=\sum_{\alpha \in \Pi_\ell}  \frac{(\alpha,\alpha )}{(\gamma, \gamma)} d_{\alpha}\alpha^\vee + \sum_{\alpha \in \Pi_{s}}d_{\alpha}\alpha^\vee.$$
Since $(\alpha,\alpha ) / (\gamma, \gamma)$ coincides with the lacing number of $\Phi$ for each long root $\alpha$,  all coefficients of $\gamma^\vee$ are multiples of the lacing number. This is a contradiction since a (co)root cannot be a multiple of an element in the (co)root lattice, since any (co)root can be a basis element of the (co)root lattice. 
\end{proof}

A direct consequence of Lemma~\ref{lem:esercizio-bou} is the following result.

\begin{pro}\label{pro:come-beta}
If all short simple roots belong to $S$, then all roots in $\Phi_{S, \be}$ share the same length. \hfill\qed \end{pro}

Lemma \ref{lem:esercizio-bou}, together with Proposition~\ref{pro:sottosistema}(\ref{item: tre}), yields a further proof of the fact that, if $\min \Phi_{S, \be}$ is short, then all roots in $\Phi_{S, \be}$ are short (see Proposition~\ref{pro:intervallo-gen}(\ref{item:3-intervallo-gen}). 
Indeed, by Proposition~\ref{pro:sottosistema}(\ref{item: tre}) $\Phi_{S, \be}$ consists of all the roots in $\Phi_{S, \ganz\be}$ that, with respect to the basis $(\Pi\smeno S)\cup \{\min \Phi_{S, \be}\}$, have the $\min \Phi_{S, \be}$-coordinate equal to $1$.

Similarly, we may obtain a criterion for $\Phi_{S, \be}$ to contain short roots in case  $\min \Phi_{S, \be}$ is long.

\begin{pro}\label{pro:min-lunga}
Let $S\cap \supp(\be)\neq \emptyset$. 
\begin{enumerate}
\item
\label{item: one} 
The set $\Phi_{S, \be}$ contains only short roots if and only if $\min \Phi_{S, \be}$ is short.
\item
\label{item: two}
The set $\Phi_{S, \be}$ contains only long roots if and only if $\min \Phi_{S, \be}$ is long and the connected component of $\min \Phi_{S, \be}$ in  the Dynkin diagram of  $(\Pi\smeno S)\cup\{\min \Phi_{S, \be}\}$ contains only long simple roots.
\item
\label{item: three}
The set $\Phi_{S, \be}$ contains both  long and short roots if and only if $\min \Phi_{S, \be}$ is long and the connected component of $\min \Phi_{S, \be}$ in  the Dynkin diagram of  $(\Pi\smeno S)\cup\{\min \Phi_{S, \be}\}$ contains some short simple root.
\end{enumerate}
\end{pro}

\begin{proof}
Part (\ref{item: one}) has already been proved.  

Parts (\ref{item: two}) and  (\ref{item: three})  follow from the following argument.

The roots in $\Phi_{S, \be}$ belong to the irreducible component of $\min \Phi_{S, \be}$ in  $\Phi_{S, \ganz \be}$. Hence,   if all simple roots in such connected component are  long, then  all roots in $\Phi_{S, \be}$ are long. On the contrary, if such connected component contains some simple short root, then the sum of all simple roots in such connected component is a root (it is  a general well-known fact that the sum of the simple roots of an irreducible root system is a root), belongs to $\Phi_{S, \be}$, and moreover, by Lemma \ref{lem:esercizio-bou},  is a short root.
\end{proof}

Lemma \ref{lem:esercizio-bou} provides a necessary condition on $\be$ in order that $\Phi_{S, \be}$ contain some long roots, namely,  that for any short simple root $\al\in S$, the coefficient $c_\al(\be)$ should be a multiple of the lacing number.
In Proposition~\ref{pro:lacing}, we see that this is also a sufficient condition, provided $S\cap \supp(\be)\neq\emptyset$. The proofs of Proposition~\ref{pro:lacing} and its  Lemma~\ref{lem:claimclaim} require the classification of the irreducible root systems. 
\par
We note that in case $S\cap \supp(\be)=\emptyset$, the above condition on $\be$ is trivially satisfied, but $\Phi_{S, \be}$ is the root subsystem generated by $\Pi\smeno S$, and hence, if $S$ contains all long simple roots, then $\Phi_{S, \be}$ only contains short roots.

\begin{lem}
\label{lem:claimclaim} 
If $\mu \in \Phi^+$ is a short root and satisfies  $\supp(\mu) = \Pi$, then there exists a short simple root $\gamma$ such that  $\mu  
+ \gamma \in \Phi^+$ and $c_{\gamma}(\mu)$ is not a multiple of the
lacing number.
\end{lem}
\begin{proof}
We prove the assertion by a  a case-by-case analysis. For each type, we list all short roots $\mu$ satisfying  $\supp(\mu) = \Pi$ and, for each such $\mu$, we provide a short simple root $\gamma$ with the required properties. The simple roots are denoted according to the Tables in \cite{Bou}.

If $\Phi$ is of type $B_l$, then $\mu = \sum_{i=1}^{l}\alpha_i$ and $\gamma=\alpha_l$.

If $\Phi$ is of type $C_l$, then $\mu = \sum_{i=1}^{j-1}\alpha_i + 2
  \sum_{i=j}^{l-1}\alpha_i  + \alpha_l$, with $1<j\leq l$, and $\gamma=\alpha_{j-1}$.

If $\Phi$ is of type $F_4$, then either $\mu = \alpha_1+\alpha_2+\alpha_3+\alpha_4$ and $\gamma= \alpha_3$, or $\mu =\alpha_1+\alpha_2+2\alpha_3+\alpha_4$ and $\gamma= \alpha_4$, or $\mu = \alpha_1+2\alpha_2+2\alpha_3+\alpha_4$ and $\gamma= \alpha_4$, or $\mu = \alpha_1+2\alpha_2+3\alpha_3+\alpha_4$ and $\gamma= \alpha_4$, or $\mu = \alpha_1+2\alpha_2+3\alpha_3+2\alpha_4$ and $\gamma= \alpha_3$.

If $\Phi$ is of type $G_2$, then either $\mu =\alpha_1+\alpha_2$, or $\mu= 2 \alpha_1+\alpha_2$, and, in both cases, $\gamma= \alpha_1$.
\end{proof}

\begin{pro}\label{pro:lacing}
Let $S\cap \supp(\be)\neq \emptyset$.
The set $\Phi_{S, \beta}$ contains long roots if and only if, for every short root $\alpha$ in $S $, the coefficient $c_\alpha(\beta)$ is a multiple of the lacing number.
\end{pro}

\begin{proof}
Let $\gamma=\sum\limits_{\alpha\in \Pi}c_\alpha(\gamma)\alpha$ be a long root in $\Phi_{S, \beta}$.
By Lemma~\ref{lem:esercizio-bou}, the coefficient  $c_\alpha(\gamma)$ is a multiple of the lacing number for every short root  $\alpha$ in $\Pi$ and so, in particular, for every short root $\alpha$ in $S$.
By the definition of $\Phi_{S, \beta}$, we have   
$c_\alpha(\beta)=c_\alpha(\gamma)$, for all $\alpha \in S$.

Conversely, suppose that $c_\alpha(\beta)$ is a multiple of the lacing number for each  short root $\alpha$  in $S$.
If $\be$ is long, then we are done. So we suppose that $\Phi$ is not simply laced and $\be$ is short. By symmetry, we may also suppose $\beta \in \Phi^+$.
Let $\mu_*$ be a maximal root  in the set $\{\mu\in\Phi_{S, \be} :  \mu \text{ is short and } \supp(\mu)\subseteq \supp(\be)\}$, which is not empty since it contains $\be$.
It is well known that the support of a root is connected and hence $\Phi\la \supp(\be) \ra $ is a standard parabolic irreducible subsystem of $\Phi$.
Thus, if $\supp(\be)$ contains only short roots, then, by the classification of the irreducible root systems, $\supp(\be)$ is a simple root system of type $A$. This is impossible since $c_{\al}(\be) >1$ for $\alpha \in S \cap \supp (\beta)$.
Hence, $\supp(\be)$ is a simple system of an irreducible not simply laced root system.
By applying Lemma \ref{lem:claimclaim} with $\mu_*$  and $\supp(\be)$ in place of $\mu$ and $\Pi$, we find a short simple root $\ga\in \supp(\be)$ such that $\mu_*+\ga$ is short and $c_\ga(\mu_*)$ is not a multiple of the lacing number.
In particular, $\ga\not\in S$ and hence $\mu_*+\ga$ still belongs to $\Phi_{S, \be}$.
By the maximality of $\mu_*$, the root $\mu_*+\ga$ is long.
\end{proof}

\section{Root systems and affine hyperplanes}
\label{sec:codim1}
In this section, we analyse  in more details   the case of codimension 1 of the  intersections  studied in Section~\ref{sec:Squalunque}. 

Let $\Phi$ be an irreducible root system (see Remark~\ref{senza loss}). Throughout this section, we fix a simple root $\alpha\in \Pi$. According to its coordinate $c_\al (\beta)$, each root  $\be=\sum\limits_{\gamma\in \Pi}c_\gamma(\beta)\gamma$ lies in one of the $2m_\al +1$ parallel affine hyperplanes $\{x\in \mathbb E : ( x , \widecheck\omega_\al) = c_\al\}$, with  $-m_\al \leq c_\al \leq m_\al$. 
To lighten the notation,  we set: 
\begin{itemize}
\item  $\Phi_{\al,k}=\Phi \cap  \{x\in \mathbb E :(  x , \widecheck\omega_\al) = k\}$ 
\item  $\Phi_{\al,\mathbb Z k}=\Phi \cap  \{x\in \mathbb E : k  \text{ divides } (  x , \widecheck\omega_\al) \}$ 
\item $\Phi_\al=\Phi \la \Pi\smeno \{\al\}\ra$
\item $W_\al=W\la\Pi\smeno \{\al\}\ra$
\item $\mu_{\al,k}=\max \Phi_{\al, k}$, for $k\in [1, m_{\alpha}]$ (see Proposition~\ref{pro:sottosistema})
\end{itemize}
In the notation of   Section~\ref{sec:Squalunque},  $\Phi_{\al,k}=\Phi_{\{\al\},\beta}$ and $\Phi_{\al,\mathbb Z k}=\Phi_{\{\al\},\mathbb Z \beta}$,
where $\beta=\sum\limits_{\gamma\in \Pi}c_\gamma(\beta)\gamma$ is any root with $c_\al (\beta)=k$. Moreover, $\Phi_\al=\Phi_{\al,0}$, $\Phi_{\al,\mathbb Z k}\supseteq \Phi_{\al,\mathbb Z k'}$ whenever $k$ divides $k'$, and $\Phi_{\al,\mathbb Z k}\cap  \Phi_{\al,\mathbb Z h} = \Phi_{\al,\mathbb Z l}$ where   $l$ is the least common multiple of $k$ and $h$. For $k=1$,  we have $\Phi_{\al, \ganz k}=\Phi$.

Notice that, by Proposition~\ref{pro:sottosistema},(\ref{item: due}), the $\Phi_{\al,k}$, with $k\neq 0$, provide a partition of $\Phi \setminus \Phi_{\alpha}$ into intervals.

By Proposition~\ref{pro:intervallo-gen} and Corollary~\ref{cor:transitivo}, we have the two following possibilities:
\begin{itemize}
\item if all roots in $\Phi_{\al,k}$ have the same length, then 
$$\Phi_{\al,k}=W_\al \cdot  \mu_{\al,k},$$
\item 
if $\Phi_{\al,k}$ contains both short and long roots, then $\mu_{\al,k}$ is long, the set $\{\be\in \Phi_{\al, k}: \be \text{ is short}\}$ has a maximum, denoted $\sigma_{\al,k}$, and  
$$\Phi_{\al,k}=W_\al\cdot  \mu_{\al,k}\, \cup\,  W_\al\cdot  \sigma_{\al,k}.$$
\end{itemize}
By the last formula in Remark \ref{rem:kostant}, however,   in both  cases
\begin{eqnarray}
\label{ugua}
\conv(\Phi_{\al,k})=\conv(W_\al\cdot  \mu_{\al,k}).
\end{eqnarray} 

Recall that   the subset $\Phi_{\al,\mathbb Z k}$, with  $k \leq m_\al$, is a root subsystem of $\Phi$, by Proposition~\ref{pro:sottosistema},(\ref{item: uno}). 
For instance, if $\Phi$ has type $G_2$, then the root subsystem $\Psi_{\al_2,k}$ has type $A_1$, $G_2$,  $A_1 \times A_1$, or  $A_2$  depending on whether $k$ has the value $0$, $1$, $2$, or $3$.

Let $\mathbb  A=\gen (\Pi\smeno \{\al\})$. 
For $x\in \mathbb E$ and $X\subseteq \mathbb E$, we let  
$\ov x\ =\pi_\mathbb  A(x)$ and $\ov X=\pi_\mathbb  A(X)$.
For all $x\in \mathbb E$, we have
$\ov x=x- \frac {(x, \widecheck\omega_\al)}{(\widecheck\omega_\al, \widecheck\omega_\al)}\widecheck\omega_\al$.
In particular,
$$\ov\ga=\ga-\frac{k}{(\widecheck\omega_\al, \widecheck\omega_\al)}\widecheck\omega_\al,$$
 for all $\ga\in \Phi_{\al, k}$.
 
Let $\ga\in \Phi_{\al, k}$ be a root of the same length as $\mu_{\al, k}$.  
The transitivity of $W_\al$ on the long roots in $\Phi_{\al, k}$ and (\ref{ugua}) imply
\begin{eqnarray}
\label{ugual}
\ov{\conv(\Phi_{\al,k})}=\conv(\ov{\Phi_{\al,k}})=\conv(W_\al\cdot \ov\ga).
\end{eqnarray}
Thus, since $\ov\ga$ is a weight of $\Phi_\al$, we have that $\ov{\conv(\Phi_{\al,k})}$ is the weight polytope of $\ov\ga$, relatively to $\Phi_\al$ in the space $\mathbb  A$. 
Clearly, $\ov\ga=0$ if and only if $\ga\in \mathbb  A^{\perp}$.

\begin{pro}\label{pro:dimensione}
Fix $\alpha \in \Pi$ and $k\in [1, m_{\alpha}]$. Let $\ga\in \Phi_{\al,k}$ be a root of the same length as $\mu_{\al, k}$. Enumerate the  irreducible components  $\Psi_1, \dots, \Psi_t$ of $\Phi_\al$ so that there exists  $h\in [1,t]$ satisfying $\ga\not\perp \Psi_i$ for $1\leq i\leq h$ and $\ga\perp \Psi_i$ for $h< i\leq t$.
Then $\ov{\Phi_{\al,k}}\subseteq \gen(\Psi_1\cup \dots\cup\Psi_h)$ and 
$$\dim(\ov{\conv(\Phi_{\al,k})})=\dim(\conv(\Phi_{\al,k}))=\mathrm{rk}(\Psi_1\cup \dots\cup \Psi_h).
$$
\end{pro}

\begin{proof}
The first equality follows from the fact that $\conv(\Phi_{\al, k})$ belongs to the affine hyperplane $\mathbb  A+\frac{k }{(\widecheck\omega_\al, \widecheck\omega_\al)}\widecheck\omega_\al$, and the restriction of $\pi_\mathbb  A$ to this affine hyperplane is the translation by $-\frac{k }{(\widecheck\omega_\al, \widecheck\omega_\al)}\widecheck\omega_\al$. 
\par
Since $\ga\perp \Psi_i$ for $h< i\leq t$, 
we have $W_\al\cdot \ga\subseteq \gen(\Psi_1\cup\dots\cup \Psi_h)+\ga$. 
By (\ref{ugual}), in order to prove the second equality, it suffices to prove that, for all $\be\in \Pi\cap(\Psi_1\cup\dots\cup\Psi_h)$, there exist $\ga_1, \ga_2\in W_\al\cdot \ga$ such that $\ga_1-\ga_2=c\be$, for some $c\neq 0$.
\par
Let $\be\in \Pi\cap \Psi_i$, for some $i\in [1, h]$. 
Since $\ga\not\perp\Psi_i$, there exist $l$ (possibly $l=1$) and a path $(\be_1, \dots, \be_l)$ in the Dynkin diagram of $ \Psi_i$ such that $\be_1=\be$, $\be_s \perp\ga$ for $s\in[1,l-1]$, and $\be_l \not\perp\ga$.  
Hence, for all $j\in [1, l-1]$, 
$$s_{\be_j}\dots s_{\be_l}(\ga)=c_j\be_j+\dots+c_l\be_l+\ga$$
for certain nonzero real $c_j, \dots, c_l$. 
We can take  the elements $s_{\be_1}\dots s_{\be_l}(\ga)$ and $s_{\be_2}\dots s_{\be_l}(\ga)$ as $\ga_1$ and $\gamma_2$, respectively.  
\end{proof}

We state the following lemma for later reference.
\begin{lem}
\label{bourba}
Let $\alpha \in \Pi$. The roots $\alpha$ and $\mu_{\al, 1}$ have the same length (the maximal length).
\end{lem}
\begin{proof}
The equality of the lengths follows from Proposition~\ref{pro:sottosistema},(\ref{item: due}), since $\al=\min\Phi_{\al, 1}$. The length is maximal by Proposition~\ref{pro:intervallo-gen},(\ref{item:3-intervallo-gen}).
\end{proof}

As a direct consequence of Proposition~\ref{pro:dimensione} and Lemma~\ref{bourba}, we obtain the following result.

\begin{cor}\label{cor:stessilegami}
Let  $\alpha \in \Pi$ and let $\Psi$ be any irreducible component of $\Phi_\al$. 
Let $k\in [1, m_{\alpha}]$ and denote by $\widetilde{ \Phi_{\al,k}}$ the set of roots of $\Phi_{\al,k}$ of the same length as $\mu_{\al, k}$.
\begin{enumerate}
\item
\label{1}  For every  $\ga, \ga'\in \widetilde{\Phi_{\al, k}}$,  we have $\ga\perp\Psi$ if and only if $\ga'\perp \Psi$.
\item
\label{2} 
For every $\ga \in \widetilde{\Phi_{\al, 1}}$, we have  $$\ga \not\perp \Psi.$$ In particular, $\dim(\conv(\Phi_{\al,1}))=\mathrm{rk}(\Phi)-1$. 
\end{enumerate}
\end{cor} 
\begin{proof}
(\ref{1}) The proof  is straightforward from Proposition~\ref{pro:dimensione}. 

(\ref{2})  Being $\Phi$  irreducible, $\alpha \not\perp \Psi$.  Since $\al$ belongs to $\Phi_{\al, 1}$ and  has the same length as $\mu_{\al, 1}$ by Lemma~\ref{bourba},  $\gamma \not\perp \Psi$ follows by (\ref{1}). The assertion on the dimension follows by Proposition~\ref{pro:dimensione}. 
\end{proof}

\section{Projections}\label{sec:inclusioni}
Let $\Phi$ be an irreducible root system. Fix an arbitrary simple root $\alpha\in \Pi$. In this section, we study the inclusion of the projected polytopes 
$\ov {\conv({\Phi_{\al, k}}})=\conv(\ov{\Phi_{\al, k}})$ in the root polytope $\mathcal P_{\Phi_\al}=\conv(\Phi_\al)$ and, more generally, the scale factors $r$ for which  $\conv(\ov{\Phi_{\al, k}})\subseteq r \cdot \mathcal P_{\Phi_\al}$. 

\begin{lem}\label{lem:-1}
Let $k\in [1, m_{\alpha}]$ and $r\in \mathbb R$.
\begin{enumerate}
\item
\label{111}
The inclusion $\conv(\ov{\Phi_{\al, k}})\subseteq r \cdot \mathcal P_{\Phi_\al}$ holds if and only if $\ov\mu\in r \cdot \mathcal P_{\Phi_\al}$ for some $\mu\in {\Phi_{\al, k}}$  of maximal length in ${\Phi_{\al, k}}$. 
\item
\label{222}
We have  $\ov{\Phi_{\al, k}}\subseteq r \cdot \mathcal P_{\Phi_\al}$ if and only if $ \ov{\Phi_{\al, -k}}\subseteq r \cdot \mathcal P_{\Phi_\al}$. 
\end{enumerate}
\end{lem}

\begin{proof}
(\ref{111}) 
By Remark~\ref{rem:kostant}, Proposition~\ref{pro:intervallo-gen} and Corollary~\ref{cor:transitivo}, the set of vertices of $\conv({\Phi_{\al, k}})$ is the set of all elements of maximal length in $\Phi_{\al, k}$, and this set is a single $W_\al$-orbit. 
Since $W_\al$ and the projection $\pi_\mathbb  A$ commute by Lemma~\ref{lem:proiezione}, and since $\Phi_\al$ is $W_\al$-stable, if a single vertex of $\conv(\ov{\Phi_{\al, k}})$ belongs to $r \cdot \mathcal P_{\Phi_\al}$, then all vertices do.
The claim follows.
\par
(\ref{222}) 
We have $\Phi_{\al, -k}=-\Phi_{\al, k}$, hence  $\ov{\Phi_{\al,-k}}=-\ov{(\Phi_{\al, k})}$, so the claim holds since $-\Phi_\al=~\Phi_\al$.
\end{proof} 

By Lemma~\ref{lem:-1},(\ref{111}), we obtain a simple criterion for checking the inclusion of $\conv(\ov{\Phi_{\al,1}})$ in the dilations of $\mathcal P_{\Phi_\al}$.

\begin{pro}\label{pro:bastaalfa}
Let $r\in \mathbb R$.
We have $\conv(\ov{\Phi_{\al, 1}})\subseteq r \cdot \mathcal P_{\Phi_\al}$ if and only if $\ov\al\in r\cdot \mathcal P_{\Phi_\al}$. 
\end{pro}

\begin{proof}
By Lemma~\ref{bourba}, the root $\al$ has maximal length in $\Phi_{\al, 1}$. The assertion follows by  Lemma~\ref{lem:-1},(\ref{111}).
\end{proof}

\begin{rem}
\label{rem:basta1}
 Let $1< k\leq m_\al$,  $\ga=\min \Phi_{\al, k}$, and $\Psi= \Phi_{\al,\mathbb Z k}$. By Proposition~\ref{pro:sottosistema}(\ref{item: two}) and~(\ref{item: tre}), the roots  $\{\ga\}\cup \Pi\smeno\{\al\}$ form a simple system for the root subsystem $\Psi$ and  $$\Phi_{\al,0}=\Psi_{\ga, 0}, \textrm{\quad while\quad }  
\Phi_{\al,k}=\Psi_{\ga, 1}.$$
The root system $\Psi$ may be reducible; however
$$
r\ov{\Psi_{\ga,1}}\subseteq \conv(\Psi_{\ga, 0})=\mathcal P_{\Psi_{\ga}}\quad\text{if and only if}\quad  r\ov{\Psi^*_{\ga,1}}\subseteq \conv(\Psi^*_{\ga, 0})=\mathcal P_{\Psi^*_{\ga}},
$$ 
where $\Psi^*$ denotes the component of $\ga$ in $\Psi$.
\end{rem}

Recall that, for each $\ve\in \Pi$, we denote by $\omega_\ve$ the fundamental weight corresponding to $\ve$, i.e. the element of $\mathbb E$ defined by the conditions $(\omega_\ve, \ve^\vee)=1$ and $(\omega_\al, \ve'^\vee)=0$ for $\ve'\in \Pi\smeno \{\ve\}$.  
Note that, for $\ve\in\Pi\smeno\{\al\}$, the projection $\ov\omega_{\ve}$ is the fundamental weight corresponding to $\ve$ in the subsystem $\Phi_\al$. 
We let $N[\alpha]$ denote the set  $\{\ve \in \Pi \smeno \{\alpha\}: \ve \not \perp \alpha \}$ of the neighbours of $\alpha$ in the Dynkin diagram of $\Phi$. 
By the following result,  $\ov{\al}$ is a negative linear combination of the  fundamental weights associated with $N[\alpha]$ in the standard parabolic subsystem $\Phi_\al$.
  
\begin{pro}
\label{pro:proiezionealfa}
The projection $\ov \alpha$ satisfies
$$\ov \alpha = \sum\limits_{\ve\in N[\al]} (\alpha, \ve^\vee ) \, \ov\omega_{\ve}.
$$
\end{pro}
\begin{proof}
The set $\{\ov\omega_\ve: \ve\in \Pi\smeno\{\al\}\}$ is the dual basis of $(\Pi\smeno\{\al\})^\vee$ in the hyperplane $\mathbb  A$. 
Hence, $\sum\limits_{\ve\in N[\al]} (\alpha, \ve^\vee ) \, \ov\omega_{\ve}$ has the same scalar products as $\al$ with the basis $(\Pi\smeno\{\al\})^\vee$ of~$\mathbb  A$. 
This implies $\sum\limits_{\ve\in N[\al]} (\alpha, \ve^\vee ) \, \ov\omega_{\ve}=\ov \al$.
\end{proof}

\begin{rem}
 If $\alpha$ is a leaf in the Dynkin diagram, the sum in Proposition~\ref{pro:proiezionealfa} consists of merely one addend, say  $(\al, \ve^\vee) \, \ov\omega_{\ve}$. Moreover, in such cases, the coefficient $ (\alpha, \ve^\vee )$  is always equal to $-1$  (and hence $\alpha$ projects to $-\ov\omega_{\ve}$) except:
\begin{itemize}
\item  in $\mathrm{C_n}$, when $\alpha=\alpha_n$, and in this case
$\ov\alpha_n= - 2 \ov\omega_{n-1}$,
\item  in $\mathrm{G_2}$, when $\alpha=\alpha_2$, and in this case 
$\ov\alpha_2= - 3 \ov\omega_1$.
\end{itemize}
Also notice that $N[\al]$ consists of at most three elements, hence $\sum\limits_{\ve\in N[\al]} (\alpha, \ve^\vee ) \, \ov\omega_{\ve}$ consists of at most three summands.
When $k=1$, we have $\al=\min \Phi_{\alpha,k}$ and $\Phi_{\alpha,\ganz k}=\Phi$. 
For $k>1$, by Proposition~\ref{pro:sottosistema}, $\Phi_{\alpha,\mathbb Zk}$ is a root system with simple system $\{\min \Phi_{\alpha,k}\}\cup \Pi\smeno\{\al\}$, hence also the projection $\ov{\min \Phi_{\alpha,k}}$ 
can be written as a linear combination of  three or less elements in $\{\ov\omega_\ve: \ve\in \Pi\smeno\{\al\}\}$. 
In fact, by a direct check, we see that, if $k>1$, then  $\ov{\min \Phi_{\alpha,k}}$
can be written as a linear combination of two or less  elements in 
$\{\ov\omega_\ve: \ve\in \Pi\smeno\{\al\}\}$. 
Otherwise, $\Pi$ and $\{\min \Phi_{\alpha,k}\}\cup \Pi\smeno\{\al\}$ would be two non isomorphic simple systems that produce isomorphic root subsystem after removing their trivalent vertex. This is impossible.
\end{rem}

\section{On a problem by Hopkins and Postnikov}
\label{sec:h-p}

In this section, we use the results in the preceding sections to tackle a problem by Hopkins and Postnikov.

In \cite{H-P}, Hopkins and Postnikov find certain applications of the root polytope of an irreducible root system $\Phi$. In particular, they need the following lemma concerning dilations of projections of the root polytope (see \cite[Lemma~2.10]{H-P}). 
\begin{lem} 
[Hopkins-Postnikov]
\label{lem:hp}
Let $U$ be a nonzero subspace of $\mathbb E$ spanned by a subset of $\Phi$. Let $\Phi_U$ be the root subsystem $\Phi \cap U$ of $\Phi$. Then there exists some $r\in \mathbb R$, with $1\leq r <2$, such that the projection on $U$ of the root polytope $\mathcal P_{\Phi}$  of $\Phi$ is contained in the dilation $r \cdot \mathcal P_{\Phi_U}$ of the root polytope of $\Phi_U$.
\end{lem}
Hopkins and Postnikov prove Lemma~\ref{lem:hp} by using Theorem~\ref{thm:faccestandard} in a  case-by-case check, and ask for a uniform proof that does not require the classification of root systems. In this section, we give a shorter proof, which unfortunately, in its last part, requires the classification (see Remark~\ref{rem: h---p}). Towards this goal, consider the following lemma.
\begin{lem}
\label{lem:uguale}
There exists some positive real number $r$, with $ r <2$, such that
$$\ov{\al}\in r \cdot \mathcal P_{\Phi_\al},$$   
for all $\alpha \in \Pi$.
\end{lem}
Lemma~\ref{lem:uguale} is in fact equivalent to Lemma~\ref{lem:hp}.
\begin{pro}
Lemma~\ref{lem:hp}  holds true for all root systems if and only if Lemma~\ref{lem:uguale} holds true for all root systems.
\end{pro}
\begin{proof}
Clearly,  Lemma~\ref{lem:hp} implies Lemma~\ref{lem:uguale}.

Conversely, suppose that Lemma~\ref{lem:uguale} holds true for all root systems. Let us prove Lemma~\ref{lem:hp}.

As explained in \cite[Appendix~A]{H-P}, we may reduce to the case when $U$ is a  maximal standard parabolic subspace, i.e., $U$ is the hyperplane generated by a codimension~1 standard parabolic subsystem of $\Phi$.  Hence the proof of Lemma~\ref{lem:hp} amounts to show that  there exists some $r\in \mathbb R$, with $1\leq r <2$, such that, for all $\alpha \in \Pi$,
\begin{eqnarray*}
\ov{\Phi_{\al, k}} \subseteq r \cdot \mathcal P_{\Phi_\al}
\end{eqnarray*}
for all $k\in [1, m_{\alpha}]$. Furthermore, Remark~\ref{rem:basta1}  implies that it suffices to show that 
there exists some $r\in \mathbb R$, with $1\leq r <2$, such that, for all $\alpha \in \Pi$,
\begin{eqnarray*}
\ov{\Phi_{\al, 1}} \subseteq r \cdot \mathcal P_{\Phi_\al}.
\end{eqnarray*}
Hence, the result follows by Proposition~\ref{pro:bastaalfa}.
\end{proof}

For each $\ve\in N[\al]$, let $\Phi^\ve$ be the component of $\ve$ in $\Phi_\al$ and $\Pi^\ve=\Pi\cap \Phi^\ve$. 
Note that,  since $\Phi$ is irreducible, for $\ve$ varying in $N[\al]$, the $\Phi^\ve$ exhaust the components of $\Phi_\al$.
Moreover, they are distinct and  at most~$3$.
Also, let $\theta^\ve$ be the highest root of $\Phi^\ve$ and, for all $\eta \in \Pi^\ve$, let $m_{\eta}$ be defined by the condition $\theta^\ve=\sum_{\eta\in \Pi^\ve} m_\eta\eta$, and let  $o_\eta=\widecheck\omega_\eta / m_\eta$.

\begin{pro}
\label{pro:somma}
Let $\Phi$ be irreducible and $\alpha \in \Pi$.
$$\min\{r: \ov{\Phi_{\al,1}}\subseteq r\cdot \mc P_{\Phi_\al}\}=\min\{r: \ov{\al}\in r\cdot \mc P_{\Phi_\al}\}=\sum\limits_{\ve\in N[\al]} - (\alpha, \ve^\vee ) c_\ve$$ 
where 
$$c_\ve = \max \{(\ov\omega_{\ve},o_{\eta}) : \textrm{$\eta$ is a $\widehat\Pi^\ve$-extremal root in $\Pi^\ve$}\}.$$
\end{pro}
\begin{proof}
The first equality follows by Proposition~\ref{pro:bastaalfa}. Let us prove the second equality.

Since $\mathcal P_{\Phi_{\al}}=-\mathcal P_{\Phi_{\al}}$, we have $\ov \al\in  r \cdot  \mathcal P_{\Phi_{\al}}$ if and only if $-\ov \al\in r \cdot \mathcal P_{\Phi_{\al}}$.
By Proposition~\ref{pro:proiezionealfa}, 
$$
-\ov \alpha =  \sum\limits_{\ve\in N[\al]} - (\alpha, \ve^\vee ) \, \ov\omega_{\ve},$$
 is dominant relatively to $\Phi_\al$ and hence, in order to check its inclusion in $ r \cdot  \mathcal P_{\Phi_{\al}} $, we may use Corollary~\ref{cor:inclusionedominanti}. If we number the elements in $N[\al]$, say $N[\al]=\{\ve_1, \dots, \ve_k\}$, by Corollary~\ref{cor:inclusionedominanti}, we obtain: 
\begin{eqnarray*}
\min\{r: \ov{\Phi_{\al,1}}\subseteq r\cdot \mc P_{\Phi_\al}\} & 
=&\max\{(-\ov\al, o_{\eta_1}+\dots+o_{\eta_k})\}\\
 & =&\max\{-(\alpha, \ve_1^\vee ) (\ov\omega_{\ve_1}, o_{\eta_1})-\dots -(\alpha, \ve_k^\vee ) (\ov\omega_{\ve_k}, o_{\eta_k})\}
 \end{eqnarray*}
where the maximum is taken over all $k$-tuples $(\eta_1, \dots, \eta_k)\in \Pi^{\ve_1}\times\dots\times \Pi^{\ve_k}$ such that  $\eta_i$ is $\widehat\Pi^{\ve_i}$-extremal. 
Hence,  the assertion follows.
\end{proof}
Note that $(\al,\ve^\vee)=-1$ except when $(\al,\al)>(\ve,\ve)$, in which case $(\al,\ve^\vee)=- (\al, \al) / (\ve, \ve)$.

\par

In order to prove Lemma~\ref{lem:uguale}, we show $\sum\limits_{\ve\in N[\al]} - (\alpha, \ve^\vee ) c_\ve <2 $, by using the classification of root systems.

Let $\Phi$ be of type  $\mathrm X_n$, with the simple roots $\{\al_1, \dots, \al_n\}$  numbered according to the Tables in \cite{Bou}. Let
$$c_i[\mathrm X_n]=\max\{(\omega_i, o_j): \textrm{$\alpha_j$ is  $\widehat\Pi$-extremal}\}.$$
By definition, $(\omega_i, o_j)=(\omega_i, \widecheck \omega_j) / m_i$. 
Furthermore,  $(\omega_i, \widecheck \omega_j)$, for $j\in[1,n]$, are the coordinates of $\omega_i$ with respect to $\{\alpha_1, \dots, \alpha_n\}$, and thus  are explicit in the tables of \cite{Bou}. Notice also that $\big( (\omega_i, \widecheck \omega_j)\big) _{i,j}$ is the inverse of the Cartan matrix $\big( (\alpha_i, \alpha_j^\vee) \big)_{i,j} $.

For $\mathrm X=\mathrm A$, all simple roots are $
\wh\Pi$-extremal and all $m_j$ are $1$, and hence $c_i[\mathrm A_n]=\max\limits_{1\leq j\leq n}\{(\omega_i, \widecheck\omega_j)\}$.
In all other cases, the $\wh\Pi$-extremal simple roots are very few. 
In the following lemma,  we list  the values of $c_i[\mathrm X_n]$ needed in the proof of Lemma~\ref{lem:uguale} ($c_i[A_n]$ is actually needed only for $i=1,2,3$).
\begin{lem}
\label{rem: computazioni}
The following equalities hold:
\begin{itemize}
\item $c_i[A_n]=i(n+1-i)/ (n+1)  $, 
\item $c_1[\mathrm B_n]=c_1[\mathrm D_n]=1$, 
\item $c_1[\mathrm C_n]=1/2$, 
\item $c_n[\mathrm B_n]=c_n[\mathrm D_n]=n/4$, 
\item $c_n[\mathrm C_n]=n/2$,
\item $c_6[\mathrm E_6]=4/3$, 
\item $c_7[\mathrm E_7]=3/2$.
\qed
\end{itemize}
\end{lem}

\begin{proof}[Proof of Lemma~\ref{lem:uguale}]
Let $\alpha \in \Pi$ and let
\begin{eqnarray}
\label{ralpha}
r_{\alpha}&=&\sum\limits_{\ve\in N[\al]} - (\alpha, \ve^\vee ) c_\ve
\end{eqnarray}
By Proposition~\ref{pro:somma},  it suffices to prove  $r_{\alpha}<2$. 
The proof consists in a case-by-case computation. 

Henceforward, $\{\al_1, \dots, \al_n\}$ is the simple system of $\Phi$ with the numbering of~\cite{Bou}, and $\al=\al_i$ for some $i\in [1,n]$.

\par
Let $\Phi\cong \mathrm A_n$. 
Then $\Phi_\al\cong \mathrm A_{i-1}\times \mathrm A_{n-i}$, where we intend $\mathrm A_{0}$ as the empty component.  
For each $\ve\in N[\al]$, we can view $\ov\omega_\ve$ as the first fundamental weight for the component of $\ve$ in $\Phi_\al$, hence $r_{\alpha}=c_1[\mathrm A_{i-1}]+c_1[\mathrm A_{n-i}]$.
\par
Let $\Phi\cong \mathrm X_n$, with $\mathrm X=\mathrm B, \mathrm C, \mathrm D$. 
Reasoning as in case $\mathrm A_n$,
 when either  $\mathrm X=\mathrm B, \mathrm C$ and $i\in[1,n-2]$, or $\mathrm X=\mathrm D$ and   $i\in [1,n-4]$, 
we find 
$$r_{\alpha}=c_1[\mathrm A_{i-1}]+c_1[\mathrm X_{n-i}].$$ 
For $\mathrm X=\mathrm B, \mathrm C$ and $i=n-1$,  we have $$r_{\alpha}=c_1[\mathrm A_{n-2}]+tc_1[\mathrm A_1]  \quad 
\text{ where }
\quad 
 t=\begin{cases} 2 & \textrm{if }X=B,\\ 
 1 & \textrm{if }X=C.
 \end{cases}
$$
For $\mathrm X=\mathrm B, \mathrm C$ and $i=n$,  we have $$r_{\alpha}=tc_1[\mathrm A_{n-1}]  \quad 
\text{ where }
\quad 
 t=\begin{cases} 1 & \textrm{if }X=B,\\ 
 2 & \textrm{if }X=C.
 \end{cases}
$$

For $\mathrm X=\mathrm D$, we have:
$$r_{\alpha}=\begin{cases} c_1[\mathrm A_{n-4}]+c_2[\mathrm A_{3}] & \textrm{if }i=n-3,\\ 
 c_1[\mathrm A_{n-3}]+c_1[\mathrm A_{1}]+c_1[\mathrm A_{1}] & \textrm{if }i=n-2,\\
 c_2[\mathrm A_{n-1}] &\textrm{if }i=n-1,n.
 \end{cases}
$$

Substituting the values of $c_i[\mathrm X]$ listed in Lemma~\ref{rem: computazioni}, we find $r_{\alpha}< 2$ in all cases.

\par

If $\Phi\cong\mathrm E_n$ ($n=6,7,8$), then
$$r_{\alpha}=\begin{cases} 
c_{n-1}[\mathrm D_{n-1}]\leq \frac{7}{4} & \textrm{if  $i=1$},\\ 
c_{3}[\mathrm A_{n-1}]\leq \frac{15}{8} & \textrm{if $i=2$},\\ 
c_{1}[\mathrm A_1]+c_{2}[\mathrm A_{n-2}]\leq
\frac{27}{14} & \textrm{if $i=3$},\\ 
c_1[\mathrm A_{1}]+c_1[\mathrm A_{2}]+c_1[\mathrm A_{n-4}]
\leq \frac{59}{30} & \textrm{if $i=4$},\\ 
c_2[\mathrm A_{4}]+c_1[\mathrm A_{n-5}]\leq 
\frac{39}{20} & \textrm{if $i=5$},\\
c_5[\mathrm D_{5}]+c_1[\mathrm A_{n-6}]\leq 
\frac{23}{12} & \textrm{if $i=6$},\\ 
c_6[\mathrm E_{6}]+c_1[\mathrm A_{n-7}]\leq 
\frac{11}{6} & \textrm{if $i=7\leq n$},\\ 
c_7[\mathrm E_{7}]=\frac{3}{2} & \textrm{if $i=8=n$}. 
\end{cases}
$$

\par
If $\Phi\cong\mathrm F_4$, then
$$r_{\alpha}=\begin{cases} 
c_3[\mathrm C_3]=\frac{3}{2} & \textrm{if  $i=1$},\\ 
c_1[\mathrm A_1]+ 2c_1[\mathrm A_2] = \frac{11}{6}& \textrm{if $i=2$},\\ 
c_1[\mathrm A_1]+c_1[\mathrm A_2]= \frac{7}{6}& \textrm{if $i=3$},\\ 
c_3[\mathrm B_3]=\frac{3}{4} & \textrm{if  $i=4$} 
\end{cases}
$$

If $\Phi\cong\mathrm G_2$, then
$$r_{\alpha}=\begin{cases} 
 c_1[\mathrm A_1]=\frac{1}{2}& \textrm{if  $i=1$},\\ 
  3c_1[\mathrm A_1]=\frac{3}{2}& \textrm{if  $i=2$}.
  \end{cases}
$$
\end{proof}

 \begin{rem}
 \label{rem: h---p}
 In  \cite{H-P}, Hopkins and Postnikov show, by a uniform argument, that Lemma~\ref{lem:hp} can be reduced
to check that, for each $\alpha \in \Pi$ and each $k\in [1,m_{\alpha}]$, the projections of the (one or two) $\Phi_{\alpha}$-dominant roots belonging $\Phi_{\alpha,k}$ fall inside  $2 \cdot   \inter{\mathcal P}_{\Phi_\al}$. This final check is done case-by-case.
Here, our case-by-case analysis amounts to just one check for each $\alpha \in \Pi$, i.e. that the uniform expression $\sum\limits_{\ve\in N[\al]} - (\alpha, \ve^\vee ) c_\ve$ of Proposition~\ref{pro:somma} specialises to a number smaller than $2$. 

 \end{rem}

\section{Enumerative formulas}
\label{sec:numerologia}
In this section, we use the results in the previous sections to give some enumerative results.

For each subset $X$ of $\Pi$, we use the shorthand notation
$W_X$ for $W\la\Pi\smeno X\ra$.
Let $\al\in \Pi$,  $W_\al=W_{\{\alpha\}}$, 
$\mathbb  A=\gen (\Pi\smeno \{\al\})$, $\pi_\mathbb  A$ be the projection on $\mathbb  A$, and, for $x\in \mathbb E$, $\ov x\ =\pi_\mathbb  A(x)$.
We recall that $N[\alpha]$ denote the set  $\{\ve \in \Pi \smeno\{\al\}: \be \not \perp \al\}$ of the neighbours of $\al$ in the Dynkin diagram of $\Pi$. 
We also set $$\widetilde N[\al]=N[\al]\cup\{\al\}.$$ 
Then $\widetilde N(\al)$ is the set of all simple roots that are not orthogonal to $\al$.

\begin{lem}
\label{m=1}
Let $\Phi$ be an irreducible root system and $\alpha \in \Pi$. The number of the roots in $\Phi_{\al,1}$ having the same length as $\al$ is 
$$\frac{|W_\al|}{|W_{\widetilde N[\al] }|}.$$
\end{lem}

\begin{proof}
By Corollary~\ref{cor:transitivo}, the parabolic subgroup $W_\al$  acts transitively both on the set  long roots of $\Phi_{\al,1}$ and on the set short roots of $\Phi_{\al,1}$. 

\par
Recall the well known fact that  the stabilizer in a Weyl group $W$ of a vector $v\in \mathcal C$ is the standard parabolic subgroup generated by the reflections through the simple roots orthogonal to $v$. 
By  Proposition~\ref{pro:proiezionealfa}, $\ov\al$ is a negative linear combination of the weigths of $\Phi_\al$ corresponding to the simple roots in $N[\al]$. 
Hence, the stabilizer in $W_\al$ of $\ov\al$, and hence of $\al$, is the parabolic subgroup $W_{\widetilde N[\al]}$. 
The assertion follows.
\end{proof}

The proof of the following result uses the classification of root systems.
\begin{pro}
\label{pro:numerodiradici}
Let $\Phi$ be a classical irreducible root system and $\al \in \Pi$ be a leaf in the Dynkin diagram.
Then, the number of the roots in $\Phi{^+}\smeno\Phi_\al$ having the same length as $\alpha$ is 
$$\frac{|W_\al|}{|W_{\widetilde N(\al)}|}.$$
\end{pro}

\begin{proof}
If $m_\al=1$, we have $\Phi^+\smeno\Phi_\al=\Phi_{\al,1}$ and we conclude by  Lemma~\ref{m=1}. 
For the classical types, this happens exactly when the leaf $\al$ is long. If $\al$ is short, we obtain the claim by duality, since $W\la  S^\vee \ra=W\la S \ra$, for each $S\subseteq \Phi^+$, and the coroots corresponding to the short roots in $\Phi^+$ whose supports contain $\alpha$ are exactly the long roots in $(\Phi^\vee)^+$ whose supports contain $\alpha^\vee$.  
\end{proof}

It is well known that the number of roots in the irreducible root system $\Phi$ equals the product of the rank of $\Phi$ times its Coxeter number (Kostant formula). 
Let $h_{\Phi}$ be the Coxeter number of $\Phi$. The following result is straightforward by Proposition~\ref{pro:numerodiradici}.

\begin{cor}
Let $\Phi$ be a root system of type  $\mathrm A_n$ or $\mathrm D_n$.  If $\alpha \in \Pi$ is a   leaf of the Dynkin diagram of $\Phi$, then
$$h_{\Phi} n - h_{\Phi_\al}(n-1)= 2 \frac{|W_\al|}{|W_{\widetilde N[\al]}|}.$$
\end{cor}

\bigskip
Let $\Phi$ be an irreducible root system and let $\{\ve_1, \ve_2, \ldots, \ve_r\}$ be the set of all simple roots of a fixed length $t$, ordered in such a way that $\Phi \la \Pi\smeno\{\ve_1, \dots, \ve_i\}\ra$ is irreducible for each $i \in [1,r]$.
Such an ordering exists, is unique in the multiple laced case, and is not unique in the simply laced case.
For each $i \in [1, r]$, let also 
$$\Pi^{(i)}=\Pi\smeno \{\ve_j:j< i\},  \qquad W^{(i)}=W\la\Pi^{(i)}\ra,$$ 
 (so $\Pi^{(1)}=\Pi$). 
By definition, $\ve_i$ is a leaf of the Dinkyn diagram of $\Pi^{(i)}$, for $i \in [1, r]$: we denote by $\ve'_i$ the unique root in $ \Pi^{(i)}$ such that 
$$\widetilde N [\ve_i]\cap \Pi^{(i)}=\{\ve_i, \ve'_i\}.$$  
Adapting the previous notation, let 
$$
W^{(i)}_{\ve_i}=W\la\Pi^{(i)}\smeno
\{\ve_i\}\ra, \qquad   W^{(i)}_{\ve_i, \ve_i'}=W\la\Pi^{(i)}\smeno
\{\ve_i, \ve'_i\}\ra.
$$

\begin{cor}
\label{iterando}
Let $\Phi$ be a classical root system, $\Phi^+_t$ be the set of all positive roots of length $t$ in $\Phi$, and  $\{\ve_1, \ve_2, \ldots, \ve_r\}$ be the set of all simple roots of length $t$, ordered in such a way that $\Phi \la\Pi\smeno\{\ve_1, \dots, \ve_i\}\ra$ is irreducible for each $i \in [1,  r]$.
Then
$$ |\Phi^+_t|=\sum\limits_{i=1}^{r} 
\frac{|W^{(i)}_{\ve_i}|}
{|W^{(i)}_{\ve_i,\ve'_i}|}.$$
\end{cor}

\begin{proof}
By definition, $\Pi^{(i+1)}=\Pi^{(i)}\smeno\{\ve_i\}$, hence, by Proposition \ref{pro:numerodiradici}, $|W^{(i)}_{\ve_i}| \, / \, |W^{(i)}_{\ve_i,\ve'_i}|$ is the number of positive roots of length $t$ in $\Phi \la \Pi^{(i)} \ra \smeno \Phi \la \Pi^{(i+1)}\ra$.
The claim follows directly, since the roots of length $t$ are contained in $\Phi\smeno \Phi\la\Pi^{(r)}\ra$.
\end{proof}

\begin{rem}
Corollary \ref{iterando} implies, in particular, that the number  $ \sum\limits_{i=1}^{r} \frac{|W^{(i)}_{\ve_i}|}
{|W^{(i)}_{\ve_i,\ve'_i}|}$ does not depend on the particular ordering of $\{\ve_1, \dots, \ve_r\}$, provided that each $\ve_i$ is a leaf of the Dynkin diagram of $\Pi\smeno\{\ve_j : j<i\}$.
This is trivial when the vertices of the  Dynkin diagram have degree   at most two. It is not trivial for type $\mathrm D$.
For type $\mathrm E$, it is false. 

Indeed, consider type $\mathrm E_6$.
If the sequence starts with the leaf (denoted $\alpha_2$ in \cite{Bou}) adjacent to the trivalent root, then  
$$ \sum\limits_{i=1}^{r} \frac{|W^{(i)}_{\ve_i}|}
{|W^{(i)}_{\ve_i,\ve'_i}|} = \frac{|W[\mathrm A_5]|}{|W[\mathrm A_2]| \times |W[\mathrm A_2]|}+ \text{ (number of positive roots in $\mathrm A_5$)}= 35,$$ 
where $W[\mathrm X_n]$ is the Weyl group of the root system of type $\mathrm X_n$.
On the other hand, if the sequence starts with another leaf,  then  
 $$ \sum\limits_{i=1}^{r} \frac{|W^{(i)}_{\ve_i}|}
{|W^{(i)}_{\ve_i,\ve'_i}|} = \frac{|W[\mathrm D_5]|}{|W[\mathrm A_4]|}+ \text{ (number of positive roots in $\mathrm D_5$)}= 36.$$
\end{rem}
 
\section*{Acknowledgments}
We thank the referee for carefully reading our manuscript and for giving helpful comments for improving the article.

\end{document}